%% file: templateArxiv.tex
\documentclass{article}
\usepackage{PRIMEarxiv}

\usepackage[utf8]{inputenc} % allow utf-8 input
\usepackage[T1]{fontenc}    % use 8-bit T1 fonts
\usepackage{url}            % simple URL typesetting
\usepackage{booktabs}       % professional-quality tables
\usepackage{nicefrac}       % compact symbols for 1/2, etc.
\usepackage{microtype}      % microtypography
\usepackage{lipsum}
\usepackage{fancyhdr}       % header
\usepackage{graphicx}       % graphics

\usepackage{caption}
\usepackage{subfig}
\usepackage{hyperref}

\usepackage{pifont}
\usepackage{amsfonts,amssymb,graphicx,listings,bm}
\usepackage{times}
\usepackage{float}
\usepackage{verbatim}
\usepackage{amscd}
\usepackage{multirow}
\usepackage{fancybox}
\usepackage{booktabs}
\usepackage{diagbox}
\usepackage{makecell}
   
\usepackage{amsmath}
\usepackage{amsthm}
\usepackage{amssymb}
\usepackage{amsfonts}
\usepackage{geometry}
\usepackage{pstricks-add}
\usepackage[all]{xy}
\usepackage{verbatim}
\usepackage{enumerate}
\usepackage{tikz}
\usepackage[ruled,vlined]{algorithm2e}

\usetikzlibrary{shapes.geometric, arrows}
\tikzstyle{arrow} = [thick,->,>=stealth]

\graphicspath{{media/}}     % organize your images and other figures under media/ folder

\input{math_commands.tex}

\input{packages}

\title{Subspace Diffusion Posterior Sampling for Travel-Time Tomography
}

\author{
  Xiang Cao$^{a}$,  Xiaoqun Zhang$^{a}$\\
  $^{a}$ School of Mathematical Sciences and Institute of Natural Sciences \\
  Shanghai Jiao Tong University \\
  Shanghai, 200240, CHINA\\
  \texttt{\url{{spawner,  xqzhang}@sjtu.edu.cn}} 
}

\begin{document}
\maketitle

\begin{abstract}
Diffusion models have been widely studied as effective generative tools for solving inverse problems. The main ideas focus on performing the reverse sampling process conditioned on noisy measurements, using well-established numerical solvers for gradient updates. Although diffusion-based sampling methods can produce high-quality reconstructions, challenges persist in nonlinear PDE-based inverse problems and sampling speed. In this work, we explore solving PDE-based travel-time tomography based on subspace diffusion generative models. Our main contributions are twofold: First, we propose a posterior sampling process for PDE-based inverse problems by solving the associated adjoint-state equation. Second, we resorted to the subspace-based dimension reduction technique for conditional sampling acceleration, enabling solving the PDE-based inverse problems from coarse to refined grids. Our numerical experiments showed satisfactory advancements in improving the travel-time imaging quality and reducing the sampling time for reconstruction.
\end{abstract}

% keywords can be removed
\keywords{Score-based diffusion models \and Diffusion posterior sampling \and  Subspace diffusion generative models \and  Nonlinear PDE-based inverse problems \and Travel-time tomography \and  Adjoint-state method; }
\section{Introduction}\label{1}
\subsection{Background}
Recently, diffusion models have emerged as powerful tools across a wide range of applications \cite{yang2023diffusion} including audio synthesis, image restoration, video generation, and so on. This success is mainly attributed to their ability to learn and leverage the implicit prior of underlying data distributions \cite{meng2021estimating}. Moreover, this capability makes diffusion models particularly suitable for solving general inverse problems.
   
In the context of solving inverse problems \cite{tarantola2005inverse, kirsch2011introduction}, we aim to recover the original parameter data $\x_0$ from the noisy or incomplete measurement $\y^{\delta}$, correlated through the forward model $\Ac$ and the measurement noise $\n$:
\begin{equation}
    \y^{\delta} = \Ac(\x_0) + \n.
\end{equation}  
Typically, the inverse problem is formulated as an optimization problem by minimizing the associated mismatch functional $\mathcal{J}(\x_0)  = \left\|\y^{\delta} -\Ac\left({\x_0}\right)\right\|^2$
under some appropriate prior information about $\x$, where $\|\cdot\|$ refers to a certain norm. In the application of diffusion models to solving inverse problems, the unconditional score function serves as this prior information to regularize the solution. Specifically, the conditional reverse sampling process for reconstruction is framed as Bayesian posterior sampling \cite{chung2023diffusion}. This approach models the conditional score function $\nabla_{\x_t} \log p_t(\x_t|\y^\delta)$ as a factorization based on Bayes' rule: 
\begin{equation}
    \nabla_{\x_t} \log p_t(\x_t|\y^\delta) = \nabla_{\x_t} \log p_t(\x_t) + \nabla_{\x_t} \log p_t(\y^\delta|\x_t)
\end{equation}
where $\nabla_{\x_t} \log p_t(\x_t)$ is approximated by utilizing the unconditional score function and $\nabla_{\x_t} \log p_t(\y^\delta|\x_t)$ is an intractable gradient term associated with the considered forward model. Various studies \cite{chung2023diffusion, song2020score, chung2022improving} have been devoted to approximating this intractable gradient term $ \nabla_{\x_t} \log p_t(\y^\delta|\x_t)$ to enable sampling from the posterior distribution $p(\x_0|\y^\delta)$. 
         
Most research in this field primarily focuses on linear/nonlinear inverse problems in non-PDE types, whose forward operators $\Ac$ can be explicitly expressed as composites that allow for automatic differentiation. For example, such linear inverse problems include inpainting \cite{lugmayr2022repaint}, deblurring \cite{chung2023diffusion}, super-resolution \cite{li2022srdiff}, Computed Tomography (CT) \cite{song2022solving}, Magnetic Resonance Imaging (MRI) \cite{chung2022score}, while nonlinear inverse problems include phase retrieval \cite{peer2023diffphase}, nonlinear deblurring \cite{chung2023diffusion} and high dynamic range \cite{mardani2023variational}. However, applying these methods directly to nonlinear PDE-based scenarios is challenging, which stems from the following two aspects:  
\begin{itemize}
    \item \textit{Nonlinear PDE-based constraints.} The inherent nonlinearity of the forward operator $\Ac$ leads to the dependence of the corresponding Fréchet derivatives $(\partial \Ac)_{\x_0}$ on the data $\x_0$, necessitating the real-time computation of the associated Jacobian matrices to evaluate the gradient of data fidelity: 
    \begin{equation}\label{constraints}
        \nabla_{\x_0} \left\|\y^{\delta} - \mathcal{A}(\x_0)\right\|_2^2 = 2(\partial \Ac)^{*}_{\x_0}\left(\mathcal{A}(\x_0) - \y^{\delta}\right)
    \end{equation}   
    under the $l_2$ norm. Given the analytical formulation of the  forward process $\Ac$, this gradient $\nabla_{\x_0} \left\|\y^{\delta} - \mathcal{A}(\x_0)\right\|_2^2$ can be efficiently evaluated using the automatic differentiation.   
    However, this approach fails due to the implicit dependence of $\Ac(\x_0)$ on $\x_0$ in the field of PDE-based inverse problems. Instead, well-established numerical solvers for both forward and adjoint PDEs are required to determine the appropriate gradient descent direction.
       
    \item \textit{Computational cost.} Conditional reverse sampling involves hundreds or even thousands of score function evaluations and gradient computations related to data fidelity. Each gradient evaluation step requires solving the forward and adjoint PDEs multiple times under different boundary conditions, making the inference much slower than the unconditional sampling. The dimensionality of both numerical PDE solvers and score functions \cite{jing2022subspace} is critical in determining the computational complexity of each sampling step, thereby influencing the overall runtime.  Furthermore, if we directly use well-established CPU-based PDE solvers \cite{jin2015finite, georgescu2013gpu, zhao2005fast, detrixhe2013parallel, wesseling1995introduction, feng2008multigrid} for convenience, the execution time for evaluating the PDE models on CPUs will dominate the sampling process, as the score model benefits from significant GPU acceleration. Therefore, accelerating the sampling process is important when using diffusion-based models to solve PDE-based inverse problems. 

\end{itemize}
Specifically, we consider solving the PDE-based travel-time tomography \cite{stefanov2019travel} problem within the diffusion posterior sampling (DPS) \cite{chung2023diffusion} framework. Here, we formulate the travel-time tomography starting from its associated forward PDE form. In the case of single point source condition, $N$ different source points $\{(u^s_n,v^s_n)\}^{N}_{n=1}$ are placed in $\Omega:= [0,1]^2$ or on $\partial \Omega$, and for each source point, the governing first travel-time field $ \Y_n(u, v)$  can be described as the following Eikonal equation \cite{rawlinson2010seismic} 
\begin{equation}\label{eikonal}
    \begin{array}{ccc}
         &\X_0(u,v) \| \nabla \Y_n(u, v)\|_2  = 1, & \forall(u, v) \in \Omega, \\
         &\text{s.t.  } \Y_n(u_n, v_n) = 0, &\text{ when } (u, v) = (u^s_n, v^s_n). 
    \end{array}
\end{equation}  
where $\|\cdot\|_2$ refers to the $L_2$ norm, $\X_0(u,v) $ refers to the positive speed field in the square $\Omega$, and its value is normalized into the range $[0,1]$. If we further denote the Eikonal operator $\mathcal{P}_n: \X_0(u,v) \longrightarrow \Y_n(u, v)$ given by \eqref{eikonal} and collect the travel-time information on the finite receiver subset $\Omega_r:= \{(u_m^r, v_m^r)\}^{M}_{m=1}\subset \partial\Omega$, the forward model for the source-receiver pair $\{(u_n^s, v_n^s), (u_m^r, v_m^r)\}$ can be stated as: 
\begin{align}\label{forward_eikonal}
    \y_{m,n}^\delta := \mathcal{P}_n(\X_0)(u_m^r, v_m^r) + \delta \bm{\eta}_{m,n},\quad  \X_0(u,v): \Rd^{2} \mapsto \Rd,
\end{align}
where $\bm{\eta}_{m,n}\in \Rd$ is an unknown noise model and $\delta$ controls the noise level. The collected data $\y_{m,n}^\delta$ from multiple source-receiver pairs constitutes a data matrix $\y^\delta \in \Rd^{M\times N}$, and the objective of travel-time tomography is to reconstruct the underlying speed field $\X_0(u,v)$ from $\y^\delta$. However, this task is a nonlinear inverse problem with severe ill-posedness. This indicates that the solution is highly sensitive to measurement noise and modeling errors, making it essential to incorporate prior knowledge of the underlying parameter distribution for more accurate and stable reconstructions.   

\begin{figure}[t!]
    \centering
    \includegraphics[width=\textwidth]{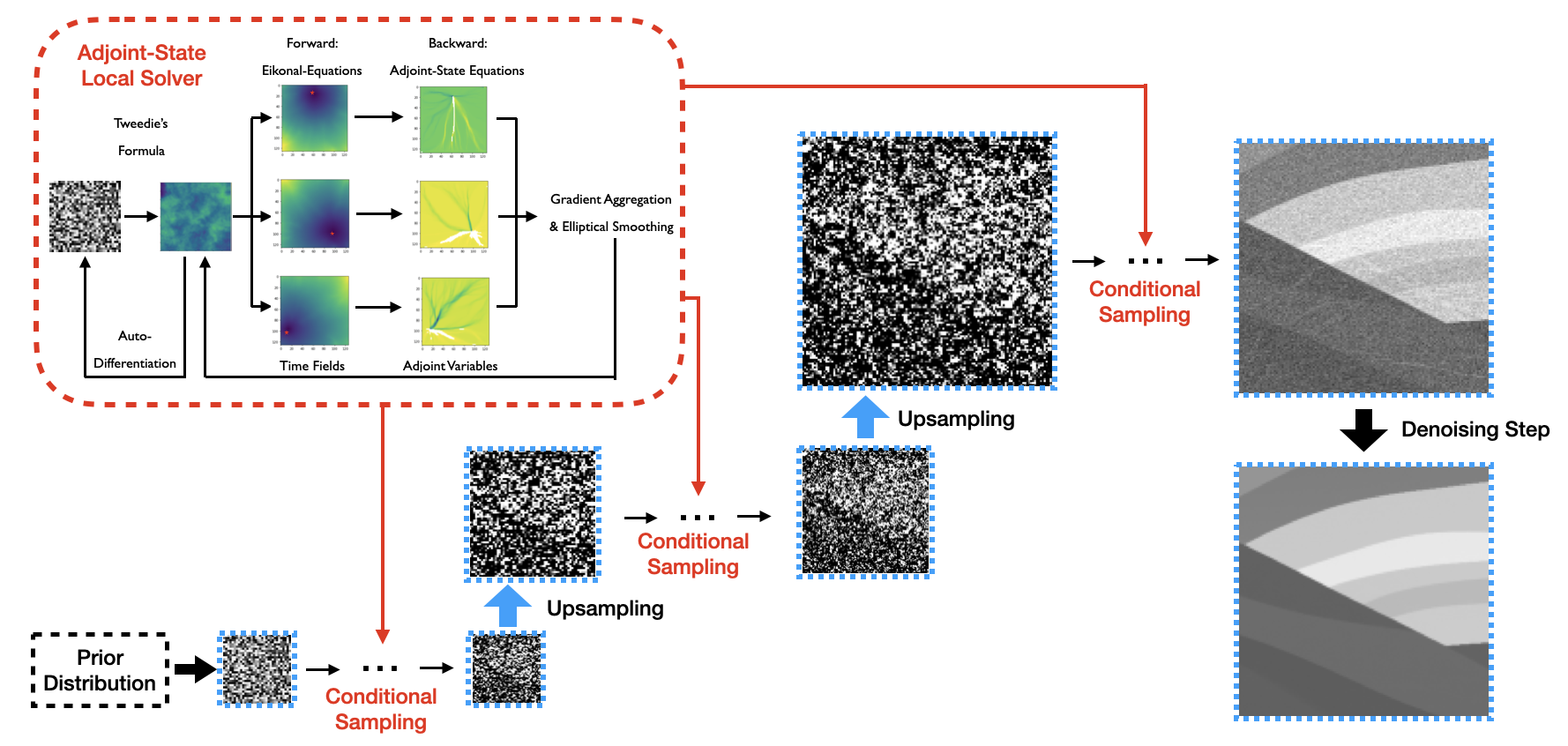}
    \caption{The Subspace-DPS framework consists of two fundamental components: the local adjoint-state solver for PDE-based gradient update and the multi-scale decomposition based on downsampling. An extra denoising step is performed for the reconstruction provided by the conditional sampling.}
    \label{workflow}
\end{figure}

\subsection{Main Contributions} 
Aiming to address the aforementioned difficulties, we develop an accelerated sampling method based on the DPS framework to solve the PDE-based travel-time tomography problem. For problem-solving, we employ the adjoint-state method with elliptical smoothing regularization \cite{leung2006adjoint} to compute the gradient of data discrepancy within the DPS framework. To accelerate sampling, we implement the subspace diffusion generative models (SGDM) \cite{jing2022subspace} for the PDE-based posterior sampling. This dimension reduction approach in the score-based framework can significantly reduce the computational dimensionality of both forward and adjoint PDE solvers, enabling solving PDE-based inverse problems from coarse to refined grids. Experimentally, we train and evaluate  score models for different resolutions on the Marmousi \cite{martin2006marmousi2} and KIT4 \cite{hamilton2018deep} datasets, each corresponding to different application scenarios of travel-time tomography. Compared to existing non-learning and learning-based methods, our proposed approach generates superior reconstruction results that better align with the prior distribution. Besides, the introduced subspace-based technique significantly accelerates the posterior sampling while maintaining comparable reconstruction accuracy.     

The remainder of this paper is structured as follows: \textbf{Section \ref{2}} provides a comprehensive literature review on solving inverse problems using score-based generative models. \textbf{Section \ref{3}} offers the brief fundamentals of diffusion posterior sampling (DPS) and subspace diffusion generative models (SDGM). Building on these techniques, \textbf{Section \ref{4}} introduces our proposed subspace posterior sampling method tailed for solving the PDE-based travel-time tomography problem. \textbf{Section \ref{5}} describes our implementation details on   
data collection and numerical settings for training and evaluation. Moreover, \textbf{Section \ref{6}} elaborates on experimental  comparison methods, and numerical results to validate our approach. Finally, \textbf{Section \ref{7}} concludes the paper with a summary of our findings and their implications for future research.

\section{Literature Review}\label{2}
Our work primarily focuses on solving inverse problems using diffusion generative models and accelerating the generation process. We will briefly review several related works on these two topics in this section.  

\textbf{Diffusion models for inverse problems.} The first category comprises methods that leverage the Tweedie formula \cite{robbins1992empirical} or other approximations for Bayesian posterior estimation. The conditional sampling process incorporates the measurement guidance through either projection \cite{kadkhodaie2021stochastic} or the gradient of data fidelity \cite{song2022solving}. The Diffusion Denoising Restoration Model (DDRM) \cite{kawar2022denoising} relies on the SVD of the forward operator to perform the conditional sampling process in the spectrum domain. The Manifold Constraint Gradient (MCG) method \cite{chung2022improving} alternates between denoising and gradient descent steps, incorporating an additional projection onto the measurement subspace to ensure data consistency. Another work $\Pi$GDM \cite{song2023pseudoinverse} introduces pseudo-inverse guidance from linear observations instead of projection. However, these methods are constrained to linear and semi-linear cases such as JPEG image restoration \cite{welker2024driftrec} and high dynamical range imaging \cite{reinhard2020high}. To overcome this limitation, the Diffusion Posterior Sampling (DPS) \cite{chung2023diffusion} method extends diffusion solvers to efficiently tackle a wide range of nonlinear inverse problems via differentiation through the score model. Nevertheless, DPS encounters difficulties in scenarios where gradients cannot be obtained through automatic differentiation, particularly in PDE-based inverse problems like \cite{10328845}. This is because the numerical operators for solving forward PDEs necessitate the use of corresponding adjoint PDE solvers for gradient back-propagation.

The second category consists of methods similar to the traditional plug-and-play prior ($\operatorname{P}^3$) approach for inverse problems \cite{venkatakrishnan2013plug}. A variational optimization method proposed by \cite{graikos2022diffusion} minimizes the data discrepancy with the diffusion error as regularization. The RED-diff approach \cite{mardani2023variational} establishes a connection between the conditional diffusion sampling and regularization-by-denoising (RED) framework \cite{romano2017little}, which allows to formulate sampling as stochastic optimization. The aforementioned two categories of methods fall under problem-agnostic diffusion models including \cite{robbins1992empirical, kawar2021snips, chung2022come, jalal2021robust, choi2021conditioning}, which can easily be adapted to different tasks without re-training the score model. Besides, another paradigm involves training a problem-specific, conditional diffusion model tailored to specific inverse problems. However, most of these models are limited to solving linear inverse problems, such as \cite{whang2022deblurring, rissanen2022generative}, due to the explicit formulation of the linear forward process.

\textbf{Sampling Acceleration.} Since hundreds or thousands of sequential evaluation steps are generally required for the reverse stochastic differential equation (RSDE) to achieve high-quality generation, score-based models have slow sampling speed. 
Considerable research \cite{bao2022analytic, dhariwal2021diffusion,
jolicoeur2021gotta, watson2021learning} have focused on reducing the number of steps in discretization and designing adaptive step-size solvers for the reverse process. Compared to sampling-based acceleration techniques, model-based acceleration approaches leverage more intricate image features or diffusion structures, thereby achieving faster inference speed while maintaining sampling accuracy. Here, we primarily enumerate four representative approaches that focus on model-based acceleration methods.
The Diffusion Probabilistic Model ODE-Solver (DPM-Solver) \cite{lu2022dpm} is designed as a fast, efficient high-order solver for equivalent probability flow ODEs to produce high-quality samples. Note that the method Denoising Diffusion Implicit Models (DDIM) \cite{song2020denoising} can be considered as a special case of DPM-Solver, which represents a more efficient class of implicit probabilistic models with non-Markovian diffusion processes. Rather than sampling in the pixel space, the Latent Diffusion Model (LDM) \cite{vahdat2021score} performs the score-based diffusion process in a low-dimensional latent space, resulting in fewer model evaluations and enhanced sample quality. Taking a different approach, the Consistency Model (CM) \cite{song2023consistency} directly maps noise to data to produce high-quality samples, with its multi-step sampling strategy further improving the quality of the generated samples. The last category consists of score-based models designed with cascaded structures, such as \cite{ho2022cascaded, saharia2022image, jing2022subspace, jiang2023low}.  Among these, the Subspace Diffusion Generative Models (SDGM) \cite{jing2022subspace} reduce computational costs for score evaluations by restricting the diffusion process through projections onto lower-dimensional subspaces. This concept shares similarities with multigrid methods \cite{wesseling1995introduction} commonly used in PDEs, where the goal is to efficiently solve problems by operating on different levels of resolution.  We draw inspiration from this approach, as lower-dimensional subspaces also imply reduced computational complexity for numerical PDE solvers, which is essential for accelerating the sampling process of diffusion posterior sampling. 

\section{Preliminaries}\label{3} 
In this section, we begin by reviewing score-based diffusion models and outlining the DPS approach for solving inverse problems within the score-based diffusion framework. Lastly, we discuss the subspace-based technique used to accelerate inference within this framework.   

\subsection{Score-based Diffusion Models}\label{3.1}

Score-based diffusion models \cite{song2020score} define the forward variance preserving stochastic differential equation (VP-SDE) for the data noising process $\x_t\in \Rd^{d},\, \forall t \in [0, T]$ as follows
\begin{equation}\label{eq:forward-sde}
    d\x_t = -\frac{\beta(t)}{2}\x_t dt + \sqrt{\beta(t)}d\w,
\end{equation}
where $\beta(t):= \beta_{\min} +  t (\beta_{\max} - \beta_{\min}) $ for $\beta_{\max},\beta_{\min}$ being two positive number serves as the noise scheduled for the forward process and $\w$ is the $d$-dimensional standard Wiener process. As the forward process progresses, the distribution of $\x_t $ ultimately converges to a standard Gaussian distribution when $t = T$. Starting from the tractable distribution $\x_T \sim \Nc(0, \mathbf{I})$, the generative process produces the data distribution $\x_0$ by solving the corresponding reverse SDE backward
\begin{equation}\label{eq:reverse-sde}
    d\x_t = \left[-\frac{\beta(t)}{2}\x_t - \beta(t)\nabla_{\x_t} \log p_t({\x_t})\right]dt + \sqrt{\beta(t)}d\bar\w,
\end{equation}
where $p_t(\x_t)$ represents the density distribution of $\x_t$ and $d\bar\w$ is the standard Wiener process running backward. Solving the reverse sampling process requires for the time-dependent score function $\nabla_{\x_t} \log p_t({\x_t})$, which is approximated by a neural network $\s_\theta(\x_t, t)$ trained with the denoising score-matching objective \cite{vincent2011connection}
\begin{equation}\label{eq:dsm}
   \min_\theta \Ed_{t} \left\{\Ed_{\x_t \sim p(\x_t|\x_0), \x_0 \sim p_{\text{data}}}\left[\|\s_\theta(\x_t, t) - \nabla_{\x_t}\log p(\x_t|\x_0)\|_2^2\right]\right\}, 
\end{equation}   
where $\x_0$ is sampled from data distribution, and the diffused sample $\x_t$ is generated by      
\begin{equation}
 \x_t = \sqrt{\bar{\alpha}(t)} \x_0 +  \sqrt{1 - \bar{\alpha}(t)} \varepsilon, \varepsilon \sim \Nc(0, \mathbf{I}), \bar{\alpha}(t) = e^{-\int_0^t\beta(s)ds}.
\end{equation}
Once the optimal parameters $\theta^*$ are obtained, various numerical solvers use $\s_{\theta^*}(\x_t, t) \simeq \nabla_{\x_t} \log p_t({\x_t})$ as a plug-in estimate to generate samples through (\ref{eq:reverse-sde}), such as \cite{song2020denoising, song2020score, ho2020denoising} and so on. 

\subsection{Diffusion Posterior Sampling for Solving Inverse Problems}\label{3.2}
Consider the general form of the forward model of an inverse problem in discrete setting expressed as
\begin{align}\label{eq:inverse}
    \y^{\delta} = \Ac(\x_0) + \bm{\n},\quad \y,\bm{\n} \in \Rd^n,\, \x \in \Rd^{d},
\end{align} 
where $\Ac(\cdot): \Rd^d \mapsto \Rd^n$ represents the forward measurement operator and $\bm{\n}$ denotes the measurement noise. The corresponding inverse problem typically suffers from severe ill-posedness and requires a robust prior for plausible reconstruction. By employing Bayes' rule, one leverages the diffusion model as prior knowledge to sample from the posterior distribution through
\begin{equation}\label{reverse-sde-posterior}
    \begin{aligned}
    d\x_t &= \left[-\frac{\beta(t)}{2}\x_t - \beta(t)(\nabla_{\x_t} \log p_t(\x_t|\y^{\delta}))\right]dt + \sqrt{\beta(t)}d\bar\w \\
    &= \left[-\frac{\beta(t)}{2}\x_t - \beta(t)(\nabla_{\x_t} \log p_t(\x_t) + \nabla_{\x_t} \log p_t(\y^{\delta}|\x_t))\right]dt + \sqrt{\beta(t)}d\bar\w. \\
    \end{aligned}
\end{equation} 
Here, the unconditional score function $\nabla_{\x_t} \log p_t(\x_t)$ serves as the prior for efficient posterior sampling in a plug-and-play fashion. In \cite{chung2023diffusion}, the intractable gradient of the log-likelihood $\nabla_{\x_t} \log p_t(\y^{\delta}|\x_t) $ is approximated through the Tweedie's approach in the following form
\begin{equation}\label{Tweedie_eq}
    \begin{aligned}
    \nabla_{\x_t}\log p(\y^{\delta}|\x_t) \simeq \nabla_{\x_t}\log p(\y^{\delta}|\hat\x_0),
    \end{aligned}
\end{equation} 
where the representation of the posterior mean $\hat\x_0$ is given by the Tweedie's formula
\begin{equation}\label{Post_eq}
    \begin{aligned}
     \hat \x_0 := \Ed[\x_0|\x_t] &= \frac{1}{\sqrt{{\bar\alpha(t)}}}(\x_t + (1 - {\bar\alpha(t)})\nabla_{\x_t} \log p_t(\x_t)).
    \end{aligned}
\end{equation}
Here, $\hat{\x}_0\left(\x_t\right)$ represents the posterior mean of $\x_0$ given  $\x_t$. We replace $\nabla_{\x_t} \log p_t(\x_t)$ with the pretrained score function $\s_{\theta^*}(\x_t,t)$ to obtain the explicit dependence of $\hat \x_0$ on $\x_t$
\begin{equation}
    \begin{aligned}
        \hat \x_0 &\simeq \frac{1}{\sqrt{{\bar\alpha(t)}}}(\x_t + (1 - {\bar\alpha(t)}) \s_{\theta^*}(\x_t,t)).
    \end{aligned}
\end{equation}
The reverse conditional sampling is now analytically tractable by substituting $\nabla_{\x_t} \log p(\x_t)$ in \eqref{reverse-sde-posterior}  with the score estimate $\s_{\theta^*}(\x_t,t)$ because the measurement distribution $p(\y^{\delta}|\x_0)$ has explicit formulation according to \eqref{eq:inverse}. For example \cite{chung2023diffusion}, in the case of the Gaussian noise with variance $\sigma^2$, the gradient of the log-likelihood term takes the form under the \(l^2\) norm
\begin{equation}\label{matrix_vector}
\begin{aligned}
    \nabla_{\x_t} \log p(\y^{\delta}|\hat\x_0(\x_t)) &= \frac{\partial \hat{\x}_0\left(\x_t\right)}{\partial \x_t} \nabla_{\hat{\x}_0} \log p(\y^{\delta}|\hat\x_0) \\
     &= - \frac{1}{\sigma^2}  \frac{\partial \hat{\x}_0\left(\x_t\right)}{\partial \x_t} 
     \nabla_{\hat{\x}_0}\left\|\y^{\delta}-\mathcal{A}\left(\hat{\x}_0\right)\right\|_2^2 \\ 
     &= - \frac{2}{\sigma^2} \frac{\partial \hat{\x}_0\left(\x_t\right)}{\partial \x_t} \left(\partial \Ac \right)^*_{\hat{\x}_0}\left(\mathcal{A}\left(\hat{\x}_0\right) - \y^{\delta}\right), 
\end{aligned}  
\end{equation}  
where $\left(\partial \Ac \right)^*_{\hat{\x}_0}$ denotes the adjoint of the Fréchet derivative at $\hat{\x}_0$. In linear scenarios, $\left(\partial \Ac \right)^*_{\hat{\x}_0}$ becomes $\mathbf{A}^T$ when $\Ac(\x_0):= \mathbf{A}\x_0$ represented as the matrix-vector product.

%where we have explicitly denoted { $\hat\x_0:=\hat\x_0(\x_t)$ to emphasize that $\hat\x_0$} is a function of $\x_t$.

\subsection{Subspace Diffusion Models}\label{Subspace Diffusion Models}
The subspace diffusion models \cite{jing2022subspace} are proposed as a dimensionality reduction technique for accelerating the sampling process. For the score-based diffusion model in the VP-SDE form \eqref{eq:forward-sde}, the corresponding subspace diffusion is defined as 
\begin{equation}\label{eq:subspace-sde}
    \begin{aligned}
         d\x_t = &\left(-\frac{\beta(t)}{2}\x_t + \sum_{k=0}^K \delta(t-t_k)(\mathbf{U}_k\mathbf{U}_k^T - \mathbf{I}_d)\x_t \right) dt \\
          + &  \sqrt{\beta(t)} \left(\sum_{k=0}^K \mathbf{1}_{[t_k,t_{k+1}]}\mathbf{U}_k\mathbf{U}_k^T \right) d\w, \\
    \end{aligned}
\end{equation}
where the time interval is divided into $0 = t_0 < t_1 < \ldots < t_K = T$, and $\mathbf{U}_k \in \mathbb{R}^{d \times d_k}$ is the column orthonormal matrix with $\mathbf{U}^T_k\mathbf{U}_k = \mathbf{I}_{d_k}$ for $d=d_0\geq d_1\geq\cdots d_{K-1}\geq d_K>0$. To constrain the diffusion process to progressively occur within lower-dimensional subspaces, one chooses $\mathbf{U}_{k}\mathbf{U}_{k}^T\mathbf{U}_{k+1}=\mathbf{U}_{k+1}$ such that the column space of $\mathbf{U}_{k+1}$ is a subspace of the column space of $\mathbf{U}_{k}$. Furthermore, denote the downsampled projection $\x^k_t:= \mathbf{U}^T_k \x_t$ and consider the diffusion process when $t\in[t_k, t_{k+1}]$
\begin{equation}\label{eq:subspace-sde-simplifed}
    \begin{aligned}
        \mathbf{U}^T_k  d\x_t = &\left(-\frac{\beta(t)}{2} \mathbf{U}^T_k\x_t +  \delta(t-t_k)\mathbf{U}^T_k (\mathbf{U}_k\mathbf{U}_k^T - \mathbf{I}_d)\x_t \right) dt \\
          + &  \sqrt{\beta(t)} \left(  \mathbf{1}_{[t_k,t_{k+1}]}\mathbf{U}^T_k\mathbf{U}_k\mathbf{U}_k^T \right) d\w. \\
    \end{aligned}
\end{equation}
Then since $d\x^k_t = \mathbf{U}^T_k  d\x_t $, the above simplifies as
\begin{equation}\label{eq:subspace-sde-simplifed1}
    \begin{aligned}
     d\x^k_t  = & -\frac{\beta(t)}{2}\x^k_t dt + \sqrt{\beta(t)}  d \w_k,
 \end{aligned}
\end{equation}
where $\w_k:= \mathbf{U}_k^T \w$ follows the standard Gaussian distribution due to the orthogonality of the $\mathbf{U}_k$ columns. Therefore, the score model $\s_k(\x^k_t, t) \simeq \nabla_{\x^k_t}\log p_t({\x^k_t})$ with lower-dimensional support can be learned using the same perturbation kernel as in the full space, thereby reducing the computational complexity of the score model evaluation.

During the sample generation process, the lower-dimensional score models $\s_k(\x^k_t, t)$ are utilized for solving the reverse SDE in the corresponding time interval
\begin{equation}\label{eq:reverse-sde-subspace}
    d\x^k_t = \left[-\frac{\beta(t)}{2}\x^k_t - \beta(t)\nabla_{\x^k_t} \log p_t({\x^k_t})\right]dt + \sqrt{\beta(t)}d\bar\w_k, \, t \in [t_k, t_{k+1}].
\end{equation}
For notation simplicity, we define the matrix notation in a similar fashion \cite{jing2022subspace}: 
\begin{itemize}
    \item $\mathbf{U}_{k\mid k-1}=  \mathbf{U}_{k-1}^T\mathbf{U}_{k} \in \mathbb{R}^{d_{k-1} \times d_{k}} $ defines the $k-1$th subspace written in the basis of the $k$th subspace.
    \item $\mathbf{P}_{k \mid k-1}=  \mathbf{U}_{k\mid k-1}\mathbf{U}^T_{k \mid k-1} \in \mathbb{R}^{d_{k-1} \times d_{k-1}} $ is the projection operator onto the $k$th subspace,
    written in the basis of the $k-1$th subspace.
\end{itemize}
To further up-sample $\x^{k-1}_{t_k}$ from $\x^{k}_{t_k}$ at the boundaries times $t_k$, we sample the orthogonal complement $\x^{k\mid k-1}_{t_k}$ from the isotropic Gaussian distribution $\mathcal{N}(0, \Sigma^{k\mid k-1}_{t_k} \mathbf{I}_{d_{k-1}})$, whose variance aims to match the marginal variance of $\x^{k\mid k-1}_{t}:= (\mathbf{I}_{d_{k-1}} - \mathbf{P}_{k \mid k-1}) \x^{k-1}_t$ at time $t_k$ 
\begin{equation} \label{eq:injection}
\begin{aligned}
    \Sigma^{k\mid k-1}_{t_k} = \frac{\bar{\alpha}(t_k)}{d_{k-1} - d_k} \mathbb{E}\left[ \lVert \x^{k\mid k-1}_{0} \rVert^2_2 \right] + 1 - \bar{\alpha}(t_k),
\end{aligned}
\end{equation}
where $\x^{k\mid k-1}_{0}$ is computed from the downsampled sample $\x^{k-1}_0$ to approximate this variance. This assumption requires the specific chosen subspaces and times such that the original magnitude of $\x^{k\mid k-1}_{0}$ is relatively small compared to the diffusion noise. In this way, we are able to sample the isotropic Gaussian noise independent of the data, and inject it via 
\begin{equation}
    \x^{k-1}_{t_k} = \mathbf{U}_{k\mid k-1} \x^{k}_{t_k} +  (\mathbf{I}_{d_{k-1}} - \mathbf{P}_{k \mid k-1}) \x^{k\mid k-1}_{t_k}.
\end{equation}
It is worth noting that while subspace diffusion shares some similarities with cascading generative models,  it preserves the flexibility of score-based models, including capabilities like conditional generation and more \cite{jing2022subspace}.    
\vspace{-0.2cm}

\section{Main Method} \label{4}

For problem-solving, we consider solving the PDE-based inverse problems within the DPS framework. Given the analytical formulation of the forward operator $\Ac$, the automatic differentiation enables efficient calculation of the gradient of the log-likelihood and the matrix-vector multiplication in \eqref{matrix_vector}.  However, in PDE-based inverse problems, \(\Ac(\hat{\x}_0)\) typically exhibits an implicit dependence on \(\hat{\x}_0\), as this relationship is governed by the underlying PDEs. Therefore, it becomes necessary to customize the gradient calculation of the data discrepancy term \(\nabla_{\hat{\x}_0}\left\|\y^{\delta}-\mathcal{A}\left(\hat{\x}_0\right)\right\|_2^2\) during the gradient back-propagation process. \textit{The adjoint-state method} will be introduced in Section \ref{Adjoint_Method} as the formulation is designed to tackle this issue.

For sampling acceleration, we propose integrating the subspace-based technique into the PDE-based diffusion posterior sampling process. Nevertheless, due to the nonlinearity of the PDE-defined forward model $\Ac$, the corresponding downsampled measurement $\y^{\delta}_k$ cannot be directly obtained from $\y^{\delta}$. The multi-resolution and discretization will be presented in Section \ref{Subspace_Downsampling} to capture the transition from the continuous multi-resolution formulation of the PDE to the discrete subspace diffusion settings. Besides, several related error bounds are proven for convergence analysis. Finally, our subspace diffusion posterior sampling strategy will be introduced in Section \ref{Subspace_Diffusion} with more details building upon the foundation laid in the previous sections. 

%As a result, the gradient \(\nabla_{\hat{\x}_0}\left\|\y^{\delta}-\mathcal{A}\left(\hat{\x}_0\right)\right\|_2^2\) cannot be directly computed via automatic differentiation, since there is no explicit analytical expression for \(\Ac(\hat{\x}_0)\).

%In Section \ref{Adjoint_Method}, we introduce the adjoint-state method for travel-time tomography, which produces the regularized gradient of the log-likelihood in posterior sampling. Section \ref{Subspace_Downsampling} presents the transition from the continuous PDE formulation to the discrete diffusion setting, defining the subspace sequences for speed fields and the corresponding discretization. In Section \ref{Subspace_Diffusion}, we describe the subspace diffusion posterior sampling approach for solving travel-time tomography, where the posterior diffusion process performs conditional sampling in each subspace. Additionally, we note that this multi-level sampling framework can be extended to other PDE-based inverse problems using the corresponding adjoint-state equations. 

\subsection{Adjoint-State Method for Travel-Time Tomography}\label{Adjoint_Method}

% Using the same notations in \eqref{eq:inverse}, we define the forward measurement operator $\Ac$ by
% \begin{equation} 
%     \Ac: \x_0 \mapsto \y_0|_{\Omega_r},
% \end{equation} 
% where $\Omega_r$ denotes the position collection of receivers on $\partial\Omega$. 

% \begin{equation}\label{J}
%      \begin{aligned}
%          \nabla_{\x_0} \mathcal{J}(\x_0) &= \frac{1}{2} \nabla_{\x_0} \left\|\y -\mathcal{A}\left(\x_0\right)\right\|_2^2 \\
%               &= (\partial{\Ac})^{*}_{\x_0}\left(\mathcal{A}\left(\x_0\right) - \y\right),
%      \end{aligned}   
% \end{equation}  
% where $\y$ is the observed noisy data on $\Omega$ and $(\partial{\Ac})^{*}_{\x_0}$ represents the adjoint operator of the linearization $\partial\Ac$ at $\x_0$. 

% Rather than utilizing the linearization approach, we adopt the adjoint-state method to compute $\nabla_{\x_0} \mathcal{J}(\x_0)$ as it offers superior computational efficiency and numerical stability. 

% This involves introducing the adjoint state variable $\boldsymbol{\lambda}$ satisfying the associated adjoint PDE
% \begin{equation}
%      \begin{aligned}
%            \mathcal{L}^*(\boldsymbol{\lambda};\y_0) = 0, \  \mathcal{L}^*(\boldsymbol{\lambda};\y_0)|_{\Omega} = \mathcal{A}\left(\x_0\right) - \y. 
%      \end{aligned}   
% \end{equation}
% By employing the adjoint-state method, we can formulate the gradient as 
% \begin{equation}
%      \begin{aligned}
%             \nabla_{\x_0} \mathcal{J}(\x_0) = \boldsymbol{\lambda}  \frac{\partial \mathcal{L}}{\partial \x_0} \,
%      \end{aligned}   
% \end{equation}        

Recall that the Eikonal-based forward model \eqref{forward_eikonal} produces the noised measurement data $\y^{\delta}_{m,n}$ for each source-receiver pair $\{(u^s_n,v^s_n),(u^r_m,v^r_m)\}$, where $(u^s_n,v^s_n)$ and $(u^r_m,v^r_m)$ represent the source location and the receiver location in $\mathbb{R}^2$, respectively. For the associated inverse problem, we aim to minimize the following misfit functional on receivers
\begin{equation}\label{misfit}
    \mathcal{J}(\X_0)=\frac{1}{2} \sum_{n=1}^N \sum_{m=1}^M \left(\y_{m,n}^\delta - \y_{m,n}\right)^2,
\end{equation}
where $\y_{m,n} = \mathcal{P}_n\left(\X_0\right)\left(u_m^r, v_m^r\right)$ denotes the travel-time obtained by the Eikonal solver $\mathcal{P}_n$ on each receiver $\left(u_m^r, v_m^r\right)$. Rather than utilizing the linearization approach, we adopt the adjoint-state method \cite{leung2006adjoint}  to compute the Fréchet derivative $\partial \mathcal{J}(\X_0)$ as it offers superior computational efficiency and numerical stability. To minimize \eqref{misfit}, we perturb the speed field $\X_0$ by $\varepsilon \tilde{\X}_0$ and the resulting change of the travel-time field $\varepsilon \tilde{\Y}_n$ satisfies the following perturbation PDE:
\begin{equation}\label{pertub}
    \nabla \Y_n \cdot \nabla \tilde{\Y}_n = - \frac{\tilde{\X}_0}{\X^3_0}.
\end{equation}
Furthermore, the change of the misfit functional is given by
\begin{equation}\label{ppp}
     \mathcal{J}(\X_0 + \varepsilon \tilde{\X}_0) -  \mathcal{J}(\X_0) = \varepsilon \sum\limits_{n=1}^N \sum\limits_{m=1}^M \left(\y_{m,n}^\delta - \y_{m,n}\right) \tilde{\y}_{m,n} +  O(\varepsilon^2), 
\end{equation}
where $\varepsilon\tilde{\y}_{m,n}$ refers to the change of the travel-time field $\varepsilon \tilde{\Y}_n$ at $(u_m^r, v_m^r)$, and the first-order term in \eqref{ppp} defines the Fréchet derivative
\begin{equation}\label{Frechet}
    \partial \mathcal{J}(\X_0)(\tilde{\X}_0) := \sum\limits_{n=1}^N \sum\limits_{m=1}^M \left(\y_{m,n}^\delta - \y_{m,n}\right) \tilde{\y}_{m,n}, 
\end{equation}
where we should choose the appropriate decent direction $\tilde{\X}_0$ such that $\partial \mathcal{J}(\X_0)(\tilde{\X}_0) < 0$. While \eqref{Frechet} presents limitations in determining this direction due to 
its implicit dependence on \(\tilde{\X}_0\),  the adjoint-state method \cite{leung2006adjoint} proves to be feasible for deriving the explicit formulation of $\partial \mathcal{J}(\X_0)(\tilde{\X}_0)$.  Here, we start by multiplying the so-called adjoint-state variable $\boldsymbol{\Lambda}_n:\Rd^2 \mapsto \Rd$ on both sides of \eqref{pertub}, and then integrate it over the whole domain $\Omega$: 
\begin{equation}\label{adjoint_edu}
\begin{aligned}
  \int_{\Omega} \frac{\tilde{\X}_0\boldsymbol{\Lambda}_n}{\X^3_0} du dv &= - \int_{\Omega}  \boldsymbol{\Lambda}_n \nabla \Y_n \cdot \nabla \tilde{\Y}_n du dv \\
   &= - \int_{\Omega}\tilde{\Y}_n  \nabla \cdot ( \boldsymbol{\Lambda}_n \nabla \Y_n) du dv + \int_{\partial \Omega}  \frac{\partial \Y_n}{\partial \boldsymbol{n}} \boldsymbol{\Lambda}_n \tilde{\Y}_n du dv, \\ 
\end{aligned}
\end{equation}
where the second equality comes from the Green formula \cite{gilbarg1977elliptic}. The adjoint-state variable $\boldsymbol{\Lambda}_n$ is introduced to eliminate the domain integral term in \eqref{adjoint_edu} and to convert the boundary integral term into $\partial \mathcal{J}(\X_0)(\tilde{\X}_0)$. To achieve this, $\boldsymbol{\Lambda}_n$ is required to satisfy the following adjoint PDE
\begin{equation}\label{adjoint_eqn} 
    \nabla \cdot (\boldsymbol{\Lambda}_n \nabla \Y_n ) = 0 
\end{equation}
with the Dirichlet boundary condition  
\begin{equation}\label{adjoint_eqn_cond} 
\frac{\partial \Y_n }{\partial \boldsymbol{n}} \boldsymbol{\Lambda}_n = \sum_{m = 1}^{M}\delta(u-u^r_m, v-v^r_m)(\y_{m,n}^\delta - \y_{m,n})
\end{equation}
where $\boldsymbol{n}$ is the unit outward normal of the boundary. We first solve the adjoint equation \eqref{adjoint_eqn} together with the boundary condition \eqref{adjoint_eqn_cond} to obtain $\Lambda_n$. Then back to \eqref{adjoint_edu}, 
\begin{equation}\label{adjoint_edu_2}
    \begin{aligned}        \int_{\Omega}\frac{\tilde{\X}_0\boldsymbol{\Lambda}_n}{\X^3_0} du dv &= - \int_{\Omega}\tilde{\Y}_n  \nabla \cdot ( \boldsymbol{\Lambda}_n \nabla \Y_n) du dv + \int_{\partial \Omega}  \frac{\partial \Y_n}{\partial \boldsymbol{n}} \boldsymbol{\Lambda}_n \tilde{\Y}_n du dv, \\ 
    &= \int_{\partial \Omega}  \left(\sum_{m = 1}^{M}\delta(u-u^r_m, v-v^r_m)(\y_{m,n}^\delta - \y_{m,n})\right) \tilde{\Y}_n  dudv \\ 
    &= \sum_{m = 1}^{M}(\y_{m,n}^\delta - \y_{m,n})\tilde{\y}_{m,n}.
    \end{aligned}
\end{equation}
Finally, we conclude that 
\begin{equation}
    \partial \mathcal{J}(\X_0)(\tilde{\X}_0) = \sum\limits_{n = 1}^{N}\int_{\Omega}\frac{\tilde{\X}_0\boldsymbol{\Lambda}_n}{\X^3_0} du dv. 
\end{equation}
However, directly using $\tilde{\X}_0 = - \sum\limits_{n=1}^{N}\boldsymbol{\Lambda}_n/{\X^3_0}$ as the update direction is inappropriate due to irregularities in the solution of \eqref{adjoint_eqn} at $(u^s_n, v^s_n)$. To obtain the well-posed descent direction, we use the elliptical smoothing regularization to choose the descent direction
\begin{equation}\label{smooth_direction}
   \tilde{\X}_0=-(\mathbf{I} -\mu \Delta)^{-1}  \left(\sum_{n=1}^N  \frac{\boldsymbol{\Lambda}_n}{\X^3_0} \right), 
\end{equation}
where $\Delta$ refers to the Laplacian operator, and the $\mu>0 $ controls the level of regularity. A homogeneous boundary condition is imposed when inverting the operator $(\mathbf{I} - \mu \Delta)$. With this specific choice of $\tilde{\X}_0$, we have
\begin{equation}
    \partial \mathcal{J}(\X_0)(\tilde{\X}_0) = - \int_{\Omega}\left( \tilde{\X}^2_0 + \mu\left| \nabla \tilde{\X}_0\right|^2\right)  du dv < 0, 
\end{equation}
which is equivalent to choosing $\tilde{\X}_0$ such that 
\begin{equation}
    \tilde{\X}_0 = \operatorname*{argmin}\limits_{\tilde{\X}}  \partial \mathcal{J}(\X_0)\left(\frac{\tilde{\X}}{\|\tilde{\X}\|_{H_\mu^1\left(\Omega\right)}}\right)
\end{equation}
in the $\mu$-weighted Sobolev space $H_\mu^1\left(\Omega\right)$.

\subsection{Multi-resolution and Discretization}\label{Subspace_Downsampling} 

Given the continuous speed field $\X_0$ as the parameter distribution of \eqref{eikonal}, we use the following notations based on convolution with downsampling and upsampling to define the multi-resolution:
\begin{equation}
    (\varphi * \X_0)(u,v) = \int_{\R^2} \varphi(\tilde{u}, \tilde{v}) \X_0(u - \tilde{u}, v - \tilde{v}) \, d\tilde{u} d\tilde{v}
\end{equation}
with downsampling and upsampling operators
\begin{equation}
     \begin{aligned}
    (\downarrow\X_0)(u,v) &:= \X_0(2u,2v),\\
    (\uparrow\X_0)(u,v) &:=  \X_0(2^{-1}u,2^{-1}v).
     \end{aligned}
\end{equation}
We continue to define the resolution reduction operator \(\Phi(\cdot)\), which maps from the \((k-1)\)-th resolution to the \(k\)-th resolution, as follows: 
\begin{equation}\label{downsample_operator} 
\X^k_0 := \Phi(\X^{k-1}_0)(u,v) = \frac{1}{2} \downarrow \left( \varphi * \X^{k-1}_0 \right), 
\end{equation} 
Similarly, the adjoint operator \(\Phi^*(\cdot)\), which maps from the \(k\)-th resolution back to the \((k-1)\)-th resolution, is defined as: \begin{equation}\label{upsample_operator} 
\Phi^*(\X^{k}_0)(u,v) := 2 \, \overline{\varphi(-\cdot, -\cdot)} * \left( \uparrow \X^{k}_0 \right),
\end{equation} 
where we denote \(\X^0_0(u,v) := \X_0(u,v)\) for \(k = 0\). Each speed field \(\X^k_0\) is extended over the entire domain \(\mathbb{R}^2\) by assuming a constant-speed background. 

In travel-time tomography, only the travel-time measurements  $\{\y^{\delta}_{m,n}\}$ are collected from a limited source-receiver pair set $\{(u_n^s, v_n^s), (u_m^r, v_m^r)\}$. To reconstruct the speed field $\X^k_0$ at the $k$-th resolution, it necessitates the corresponding measurements $\y^{\delta}_{k,m,n}$ at the downsampled receiver location $(\frac{1}{2^k}u_m^r, \frac{1}{2^k}v_m^r)$ for reconstruction. However, the associated travel-time field \(\Y^k_n\) satisfies the nonlinear Eikonal equation:
\begin{equation}\label{eikonal_n}
\X^k_0(u, v) \, \| \nabla \Y^k_n(u, v) \|_2 = 1, 
\end{equation}
which indicates that $\y^{\delta}_{m,n}$ can only serve as an estimation of $\y^{\delta}_{k,m,n}$ based on the following approximation: 
\[ 
\begin{aligned} \y_{k,m,n} &= \mathcal{P}_n(\X^k_0)\left(\frac{1}{2^k}u_m^r, \frac{1}{2^k}v_m^r\right), \\ & \approx \mathcal{P}_n(\X_0)(u_m^r, v_m^r), \\ & = \y_{m,n}. \end{aligned} 
\] 
Thus, for the corresponding inverse problem, we aim to minimize the surrogate misfit functional:
\begin{equation}\label{surrogate_misfit_k}
    \mathcal{J}_k(\X^k_0) = \frac{1}{2} \sum_{n=1}^N \sum_{m=1}^M \left(\y^{\delta}_{m,n} - \y_{k,m,n}\right)^2
\end{equation}
by replacing $\y^{\delta}_{k,m,n}$ with $\y^{\delta}_{m,n}$. Additionally, we can derive a closed-form upper bound for the gap between $\mathcal{J}_k(\X^k_0)$ and $\mathcal{J}(\X_0)$, as presented in the following theorem:
\begin{restatable}[]{theorem}{bound}\label{bound}   
For the given convolution kernel $\varphi(u,v)$ with $\left|\prod\limits_{l=0}^{k-1}\hat{\varphi}(2^{l}\xi, 2^{l}\eta) - 1 \right| \le \varepsilon \left(1 + \mu(|\xi|^2 + |\eta|^2)\right)^{\frac{r}{2}}$, the gap between  $\mathcal{J}_k(\X^k_0)$ and  $\mathcal{J}(\X_0)$ is upper bounded by
\begin{equation}\label{bound_eq}
    |\mathcal{J}_k(\X^k_0) - \mathcal{J}(\X_0)| \le  \varepsilon \|\tilde{\X}_0\|^2_{H_\mu^1\left(\Omega\right)} \|\X_0\|_{H_{\mu}^{1+r}\left(\Omega\right)} + O(\varepsilon^2)
\end{equation}
where $\widehat{\varphi}(\cdot)$ denotes the 2D Fourier transform of $\varphi(\cdot)$, and $\tilde{\X}_0$ refers to the elliptical smoothed perturbation given by the adjoint-state method in \eqref{smooth_direction}.
\end{restatable}

Since the diffusion generative process is defined on the discrete grids, we should consider the discretization of the $k$-th resolution speed field samples $\X^k_0$ as the projection onto the subspaces. For the mesh size $ h = 2^{-s}$, let us consider the regular grid 
\begin{equation}\label{regular_grid}
    \mathbb{Z}^2_{h} = \{(i h, j h)|i,j\in \mathbb{Z}\}.  
\end{equation}
Next, denote $\Omega^k_{h} = \mathbb{Z}^2_{h} \cap [0, 2^{-k})^2$ to be the inner grid points in $[0, 2^{-k})^2 $, and obtain the corresponding discretization $\x_0^k\in \R^{2^{s-k}\times 2^{s-k}}$ via
\begin{equation}
    \x_0^k[i,j] = \X^k_0(i h, j h), \text{ } \forall (i h, j h) \in  \Omega^k_{ h}.
\end{equation}
Then, we can directly formulate $\x_0^{k+1}$ from $\x_0^{k}$ by the following theorem:
\begin{restatable}[]{theorem}{kernel}\label{kernel}
For the given convolution kernel $\varphi_{h}(u,v) = \frac{1}{4}\sum\limits_{(a,b)\in \{0,1\}^2} \delta(u+ah,v+bh)$, the following relationship holds between $\x_0^k$ and $\x_0^{k-1}$:
\begin{equation}\label{kernel_eq}
    \x_0^{k+1}[i,j] =  \frac{1}{2}\mathcal{D}_{k}\x_0^{k}[i,j] = \frac{1}{2}\left(\sum\limits_{(a,b)\in \{0,1\}^2} \frac{\x_0^{k}[2i+a, 2j+b]}{4}\right)
\end{equation}
where $\mathcal{D}_k: \mathbb{R}^{2^{s-k} \times 2^{s-k}} \rightarrow \mathbb{R}^{2^{s-k-1} \times 2^{s-k-1}}$ represents the average-pooling operator.
\end{restatable}
\textbf{Remarks.} In the subspace diffusion model, the row space of the orthogonal matrix $\mathbf{U}^T_k$ defines the projection subspace. Here, we define the orthogonal matrix $\mathbf{U}^T_k$ via
\begin{equation} \label{eq:uk}
    \mathbf{U}_k^T \x_0 := \x^k_0 = \frac{1}{2^k} \left(\prod\limits_{l=0}^{k-1} \mathcal{D}_l\right) \x_0
\end{equation}
where $\x_0$ is represented as the column vector and the matrix in $\mathbf{U}_k^T \x_0$ and $\{\mathcal{D}_l\x^l_0\}_{l=0}^{k-1}$, respectively. With this specific choice of $\varphi_{h}(u,v)$, we can formulate the upper bound of the gap between $\mathcal{J}_k(\X^k_0)$ and  $\mathcal{J}(\X_0)$ as follows:
\begin{restatable}[]{corollary}{bounddiscrete}\label{bound_discrete}
Given the convolution kernel $\varphi_{h}(u,v)$ with enough small $h$ in Theorem \ref{kernel}, the gap between $\mathcal{J}_k(\X^k_0)$ and  $\mathcal{J}(\X_0)$ is upper bounded by
\begin{equation}\label{bound_discrete_eq}
    |\mathcal{J}_k(\X^k_0) - \mathcal{J}(\X_0)| \le C_{k,\mu} h \|\tilde{\X}_0\|^2_{H_\mu^1\left(\Omega\right)} \|\X_0\|_{H_{\mu}^{2}\left(\Omega\right)} + O(h^2)
\end{equation}
where $C_{k,\mu}$ is a fixed constant given $k, \mu$, and $\tilde{\X}_0$ refers to the elliptical smoothed perturbation given by the adjoint-state method in \eqref{smooth_direction}.
\end{restatable}
Finally, let us denote four neighbors of $[i,j]$ by $\mathcal{N}_{[i,j]}$, and we can further approximate the numerical solution of $\Y_n^k$ in \eqref{eikonal_n} using a discretized array $\y^k_n \in \R^{2^{s-k}\times 2^{s-k}}$ as follows:
\begin{equation}\label{local_solver_section}
   \sum\limits_{[a,b]\in \mathcal{N}_{[i,j]}} \left[\left(\frac{\y^k_n[i,j] - \y^k_n[a,b]}{h}\right)^{+}\right]^2 = \left(\frac{1}{\x^k_0[i,j]}\right)^2,
\end{equation}
\begin{equation}\label{start_condition_section}
    \text{s.t. } \y^k_n\left[i_n^{s,k}, j_n^{s,k}\right] = 0,
\end{equation}
where $\left[i_n^{s,k}, j_n^{s,k}\right] = \left[\lfloor \frac{u_n^s}{2^k h} \rfloor,\lfloor\frac{v_n^s}{2^k h}\rfloor\right]$ denotes the source grid point, and $(x)^+$ is the positive part of $x$. Note that this local solver with the upwind scheme is designed to update the neighbors of $[i,j]$, and we typically use the fast marching method (FMM) \cite{sethian1999fast} to solve \eqref{local_solver_section} \eqref{start_condition_section} over the whole domain $\Omega^k_h$. To numerically calculate the adjoint-state variable $\boldsymbol{\Lambda}^k_n$ in \eqref{adjoint_eqn} \eqref{adjoint_eqn_cond}, we use the discretization $\boldsymbol{\lambda}_n^k$ for approximation 
\begin{equation}\label{adjoint_discrete}
  \sum\limits_{[a,b]\in \mathcal{N}_{[i,j]}}\left[ \left(\frac{\y^k_n[i,j] - \y^k_n[a,b]}{h}\right)^{+} \frac{\boldsymbol{\lambda}^k_n[i,j]}{h} - \left(\frac{\y^k_n[a,b] -\y^k_n[i,j]}{h}\right)^{+}\frac{\boldsymbol{\lambda}^k_n[a,b]}{h}\right] = 0,
\end{equation}
\begin{equation}\label{adjoint_discrete_cond}
     \text{s.t. } \left(\frac{\y^k_n[p,q] - \y^k_n[i^{r,k}_m,j^{r,k}_m]}{h}\right) \boldsymbol{\lambda}_n[i^{r,k}_m,j^{r,k}_m]  = \y^{\delta}_{m,n} - \y_n[i^{r,k}_m,j^{r,k}_m],
\end{equation}
where $\left[i^{r,k}_m, j^{r,k}_m\right] = \left[\lfloor \frac{u^{r}_m}{2^kh} \rfloor,\lfloor\frac{v^{r}_m}{2^kh}\rfloor\right]$ denotes the receiver grid point, and $[p,q] \in \mathcal{N}_{[i,j]}$ is off the boundary. Note that the discrete adjoint equation \eqref{adjoint_discrete} \eqref{adjoint_discrete_cond} can be solved at each point following the reverse pop order provided by FMM \cite{deckelnick2011numerical}.

For the corresponding inverse problem, the misfit function $\mathcal{J}^h_k(\cdot): \R^{2^{s-k}\times 2^{s-k}} \mapsto \R$: 
\begin{equation}\label{misfit_approx_discrete}
    \mathcal{J}^h_k(\x^k_0)=\frac{1}{2} \sum_{n=1}^N \sum_{m=1}^M \left(\y^{\delta}_{m,n} - \y^k_n[i^{r,k}_m,j^{r,k}_m]\right)^2
\end{equation}
is defined on the discrete sample $\x_0^k$. With the help of Theorem 2.6 in \cite{deckelnick2011numerical}, we can derive the upper bound estimate of the gap between $\mathcal{J}^h_k(\x^k_0)$ and  $\mathcal{J}(\X_0)$ by combining the result provided by Corollary \ref{bound_discrete}:
\begin{restatable}[]{theorem}{bounded}\label{bound_ed}
For the given convolution kernel $\varphi_{h}(u,v)$ in Theorem \ref{kernel}, the gap between $\mathcal{J}^h_k(\x^k_0)$ and  $\mathcal{J}(\X_0)$ is upper bounded by
\begin{equation}\label{bound_ed_eq}
    |\mathcal{J}^h_k(\x^k_0) - \mathcal{J}(\X_0)| \le  C_{k,\mu} h \|\tilde{\X}_0\|^2_{H_\mu^1\left(\Omega\right)} \|\X_0\|_{H_{\mu}^{2}\left(\Omega\right)} + C_1 \sqrt{h} + O(h)
\end{equation}
    where $C_{k,\mu}$ is a fixed constant given $k,\mu$, and $C_1$ depends on $\X_0^k$ and $\y^{\delta}$. 
\end{restatable}
\textbf{Remarks.} As $h \rightarrow 0$, Theorem \ref{bound_ed} guarantees the numerical consistency for optimization, that is, the surrogate misfit function satisfies $\lim\limits_{h \rightarrow 0} \mathcal{J}^h_k(\x^k_0) = \mathcal{J}(\X_0)$.

% To numerically calculate the adjoint-state variable $\boldsymbol{\Lambda}_n^k$ in \eqref{adjoint_eqn} \eqref{adjoint_eqn_cond}, we use the discretization $\boldsymbol{\lambda}_n^k$ for approximation 
% \begin{equation}\label{adjoint_discrete}
%   \sum\limits_{[i,j]\in \mathcal{N}_{[a,b]}}\left[ \left(\frac{\y^k_n[a,b] - \y^k_n[i,j]}{h}\right)^{+} \frac{\boldsymbol{\lambda}_n^k[a,b]}{h} - \left(\frac{\y^k_n[i,j] -\y^k_n[a,b]}{h}\right)^{+}\frac{\boldsymbol{\lambda}_n^k[i,j]}{h}\right] = 0,
% \end{equation}
% \begin{equation}\label{adjoint_discrete_cond}
%      \text{s.t. } \left(\frac{\y^k_n[p,q] - \y^k_n[a^r_m,b^r_m]}{h}\right) \boldsymbol{\lambda}_n^k[a^r_m,b^r_m]  = \y^{\delta}_{m,n} - \y^k_n[a^r_m,b^r_m],
% \end{equation}
% where $\left[a^r_m, b^r_m\right] = \left[\lfloor \frac{u^r_m}{2^k h} \rfloor,\lfloor\frac{v^r_m}{2^k h}\rfloor\right]$ denotes the receiver index on $\x^k_0$, and $[p,q] \in \mathcal{N}_{[a,b]}$ is off the boundary. Note that the discrete adjoint equation \eqref{adjoint_discrete} \eqref{adjoint_discrete_cond} can be solved at each point following the reverse pop order provided by FMM \cite{deckelnick2011numerical}.

\subsection{{Subspace Diffusion Posterior Sampling}}\label{Subspace_Diffusion}

% Generally, the gradient of the log-likelihood term in \eqref{reverse-sde-posterior} is approximated by 
% \begin{equation}\label{matrix_vector}
% \begin{aligned}
%     \nabla_{\x_t} \log p(\y^{\delta}|\hat\x_0(\x_t)) &= \frac{\partial \hat{\x}_0\left(\x_t\right)}{\partial \x_t} \nabla_{\hat{\x}_0} \log p(\y^{\delta}|\hat\x_0) \\
%      &= -\rho  \frac{\partial \hat{\x}_0\left(\x_t\right)}{\partial \x_t} 
%      \nabla_{\hat{\x}_0}\left\|\y^{\delta}-\mathcal{A}\left(\hat{\x}_0\right)\right\|_2^2,
% \end{aligned}  
% \end{equation}  
% where the log-likelihood \(\log p(\y^{\delta}|\hat{\x}_0)\) is expressed as \(-\rho \|\y^{\delta} - \Ac(\hat{\x}_0(\x_t))\|_2^2\) under the \(l^2\) norm in the case of Gaussian noise. 

For the travel-time tomography problem, the forward operator $\mathcal{A}$ is defined by
\begin{equation}
    \mathcal{A}: \x_0 \in \mathbb{R}^{2^s \times 2^s} \mapsto \y := \left[\y_n\left[i_m^r, j_m^r\right]\right]_{m=1,n=1}^{M,N} \in \mathbb{R}^{M \times N}, 
\end{equation}
where $\y_n$ refers to the discrete travel-time field solved by \eqref{local_solver_section} when $k=0$, and $(i_m^r, j_m^r)$ represents the receiver grid point as indicted in \eqref{misfit_approx_discrete} when $k=0$. Here, we consider utilizing the discrete adjoint-state method to provide $\nabla_{\hat{\x}_0} \left\|\y^{\delta}-\mathcal{A}\left(\hat{\x}_0\right)\right\|_2^2$ for guidance during the diffusion posterior sampling process, that is
\begin{equation}\label{dps_evaluation}
\begin{aligned}
    \nabla_{\x_t} \log p(\y^{\delta}|\hat\x_0(\x_t)) &= - 2\rho \frac{\partial \hat{\x}_0\left(\x_t\right)}{\partial \x_t} \nabla_{\hat{\x}_0}\mathcal{J}^h_0(\hat{\x}_0). \\
    & = 2 \rho \frac{\partial \hat{\x}_0\left(\x_t\right)}{\partial \x_t} \left(\sum\limits_{n = 1}^{N}\boldsymbol{\lambda}_n /\hat{\x}^3_0\right).
\end{aligned}  
\end{equation}
Note that $\mathcal{J}^h_0(\hat{\x}_0)$ is the misfit function defined in \eqref{misfit_approx_discrete}, and the discrete adjoint-state variable $\boldsymbol{\lambda}_n$ is solved in \eqref{adjoint_discrete} when $k = 0$. As suggested in the adjoint-state method, $\boldsymbol{\lambda}_n$ exhibits singularities at the source points, leading the numerical instabilities for the matrix-vector multiplication $\frac{\partial \hat{\x}_0\left(\x_t\right)}{\partial \x_t} \nabla_{\hat{\x}_0}\mathcal{J}^h_0(\hat{\x}_0)$. To overcome this, we resort to the elliptical smoothing regularization to perform the gradient via
\begin{equation}\label{dps_evaluation_smoothed}
\begin{aligned}
    \nabla_{\x_t} \log p(\y^{\delta}|\hat\x_0(\x_t)) &\simeq - 2\rho \frac{\partial \hat{\x}_0\left(\x_t\right)}{\partial \x_t} (\mathbf{I} -\mu \Delta)^{-1}\left(\nabla_{\hat{\x}_0}\mathcal{J}^h_0(\hat{\x}_0)\right) \\
    & = 2 \rho \frac{\partial \hat{\x}_0\left(\x_t\right)}{\partial \x_t}  (\mathbf{I} -\mu \Delta)^{-1} \left(\sum\limits_{n = 1}^{N}\boldsymbol{\lambda}_n /\hat{\x}^3_0\right)
\end{aligned}  
\end{equation}
where $(\mathbf{I} -\mu \Delta)^{-1}$ is defined in \eqref{smooth_direction} with a homogeneous boundary condition. Note that utilizing the adjoint-state method to compute the gradient guidance indeed serves as a plug-and-play approach for conditional sampling. To accelerate the sampling process, we consider the reverse posterior SDE with the same form 
\begin{equation}\label{eq:reverse-subspace-sde-posterior} 
    \begin{aligned}
     d\x^k_t = \left[-\frac{\beta(t)}{2}\x^k_t - \beta(t) \left( \nabla_{\x^k_t} \log p_t({\x^k_t}) +  \nabla_{\x^k_t} \log p(\y^\delta|\hat\x^k_0(\x^k_t))  \right)\right]dt + \sqrt{\beta(t)}d\bar\w_k, 
    \end{aligned}
\end{equation}
for $t \in [t_k, t_{k-1}]$ in each subspace. Besides, $\x^k_t$ is obtained via 
\begin{equation} \label{eq:uk_here}
     \x^k_t := \mathbf{U}_k^T\x_t  = \frac{1}{2^k} \left(\prod\limits_{l=0}^{k-1} \mathcal{D}_l\right)\x_t
\end{equation}  
where the average-pooling operator $\mathcal{D}_k$ is defined in \eqref{kernel} when the convolution kernel $\varphi_h$ is given as the condition in Theorem \ref{kernel}. As indicated in the aforementioned subsection, the true downsampled observation $\y^\delta_{k,m,n}$ is unobtainable due to the nonlinearity of the Eikonal equation. Therefore, we choose the surrogate misfit function $\mathcal{J}^h_k(\cdot)$ for sampling, that is
\begin{equation}\label{dps_k}
    \nabla_{\x^k_t} \log p(\y^{\delta}|\hat\x^k_0(\x^k_t)) \simeq  - 2\rho_k \frac{\partial \hat{\x}^k_0\left(\x^k_t\right)}{\partial \x^k_t} (\mathbf{I} -\mu \Delta)^{-1}\left(\nabla_{\hat{\x}^k_0}\mathcal{J}^h_k(\hat{\x}^k_0)\right),
\end{equation}
where $\mathcal{J}^h_k(\cdot)$ is defined in \eqref{misfit_approx_discrete}, and Theorem \ref{bound_ed} guarantees that $\mathcal{J}^h_k(\cdot)$ approximates the true misfit function $\mathcal{J}(\cdot)$ as $h$ is small enough. Besides, according to the subspace diffusion process in \eqref{eq:subspace-sde-simplifed1}, we can treat the forward diffusion process as occurring independently for each downsampled sample $\x^k_0$. Then, the subspace score model $\s_k(\x^k_t, t)$ corresponds to the score model trained on a lower-resolution variant of the original dataset. The downsampled posterior mean $\hat{\x}^k_0\left(\x^k_t\right)$ is then approximated via the Tweedie formula as well
\begin{equation}
    \begin{aligned}
        \hat \x^k_0 &= \frac{1}{\sqrt{{\bar\alpha(t)}}}(\x^k_t + (1 - {\bar\alpha(t)})\nabla_{\x^k_t} \log p_t(\x^k_t)) \\
        &\simeq \frac{1}{\sqrt{{\bar\alpha(t)}}}(\x^k_t + (1 - {\bar\alpha(t)}) \s_k(\x^k_t, t)),
    \end{aligned}
\end{equation}
Besides, the up-sampling procedure at the boundary times $t_k$ follows the same approach as the subspace diffusion models introduced in Section \ref{Subspace Diffusion Models}.

We should note that the primary benefit of the subspace-based diffusion posterior sampling lies in its capability to simultaneously reduce the dimensionality of score function evaluations and the complexity of solving numerical PDEs. For example, in the travel-time tomography, the computational complexity of the fast marching method \cite{sethian1999fast} and the fast sweeping method \cite{zhao2005fast} for solving the Eikonal equations requires \(O(d \log d)\) and \(O(d)\), respectively. Although the computational complexity of evaluating the score model is $O(d)$, the GPU can significantly accelerate the real-time computation of the score function. If we directly plug well-established CPU-based solvers in for convenience, consequently, the time spent on evaluating the PDE-based gradient becomes the dominant factor in each sampling step. This can be a bottleneck that limits the posterior sampling speed if no acceleration techniques are applied. Furthermore, it is important to emphasize that this subspace posterior sampling approach is \textbf{applicable to other PDE-based inverse problems}, as long as the adjoint equations can be explicitly formulated. Here, we summarize the overall algorithm in Algorithm \ref{alg:sampling}.

 \begin{algorithm}[t]
\caption{Conditional subspace sampling for solving the travel-time tomography problem}\label{alg:sampling}
\KwIn{projection matrices $\mathbf{U}^T_k$, projection times $t_k$, observed travel-time data $\y^\delta \in \mathbb{R}^{M\times N}$, the Eikonal-based forward operator $\Ac$, score models $\s_k(\x^k_t,t)$, $k = 0,\ldots, K-1$}  
\KwOut{reconstruction sample $\x_0$ from $p(\x_0|\y^\delta)$} 
$\x^{K-1}_T$ $\gets$ sample from $\mathcal{N}(0, \mathbf{I}) \in \mathbb{R}^{2^{s-K+1} \times 2^{s-K+1}}$\;
 
\For{$k\leftarrow K-1$ \KwTo $0$}{
  \For{ $t \in [t_{k+1},t_{k}]$}{$\hat{\x}^k_0 \gets \frac{1}{\sqrt{\bar\alpha(t)}}(\x^{k}_{t} + {(1 - \bar\alpha(t))}\s_k(\x^{k}_{t}, t))$ \;
  
 {$(\mathbf{I} -\mu \Delta)^{-1}\left(\nabla_{\hat{\x}^k_0}\mathcal{J}^h_k(\hat{\x}^k_0)\right)$}  $\gets$
  compute the regularized gradient using the adjoint-state method in \eqref{dps_evaluation_smoothed}\;
 {$\nabla_{\x^k_t} \log p(\y_k|\hat\x^k_0(\x^k_t))$} $\gets$ compute the matrix-vector multiplication in \eqref{dps_k} via auto-differentiation \;
  }   

 $\x^{k}_{t_k} \gets$ solve the conditional reverse sampling in \eqref{eq:reverse-subspace-sde-posterior} from $\x^{k}_{t_{k+1}}$ in $ [t_{k+1},t_{k}]$\;
  
   \If{$k > 0$}{
    $ \x^{k\mid k-1}_{t_k} \gets $ sample from $\mathcal{N}(\mathbf{0}, \Sigma^{k\mid k-1}_{t_k} \mathbf{I}) \in \mathbb{R}^{2^{s-k} \times 2^{s-k}}$  \;
    
    $\x^{k\mid k-1}_{t_k} \gets (\mathbf{I} - \mathbf{P}_{k \mid k-1})\x^{k\mid k-1}_{t_k} $ \;
    
    $\x^{k-1}_{t_k} \gets \mathbf{U}_{k\mid k-1}\x^{k}_{t_k} +\x^{k\mid k-1}_{t_k}$ \;
    }
  }
$\x_{0}$ $\gets$ denoise $\x^0_0$ via $\frac{1}{\sqrt{\bar\alpha(\varepsilon)}}(\x^0_0 + {(1 -\bar\alpha(\varepsilon))}\s_0(\x^0_{0}, \varepsilon))$ where $\frac{\varepsilon}{T}= 10^{-3}$\; 
\end{algorithm}

\section{Implementation details}\label{5}
\subsection{Dataset Collection}  
\begin{figure}[t!]
    \centering
    \includegraphics[width=0.9\textwidth]{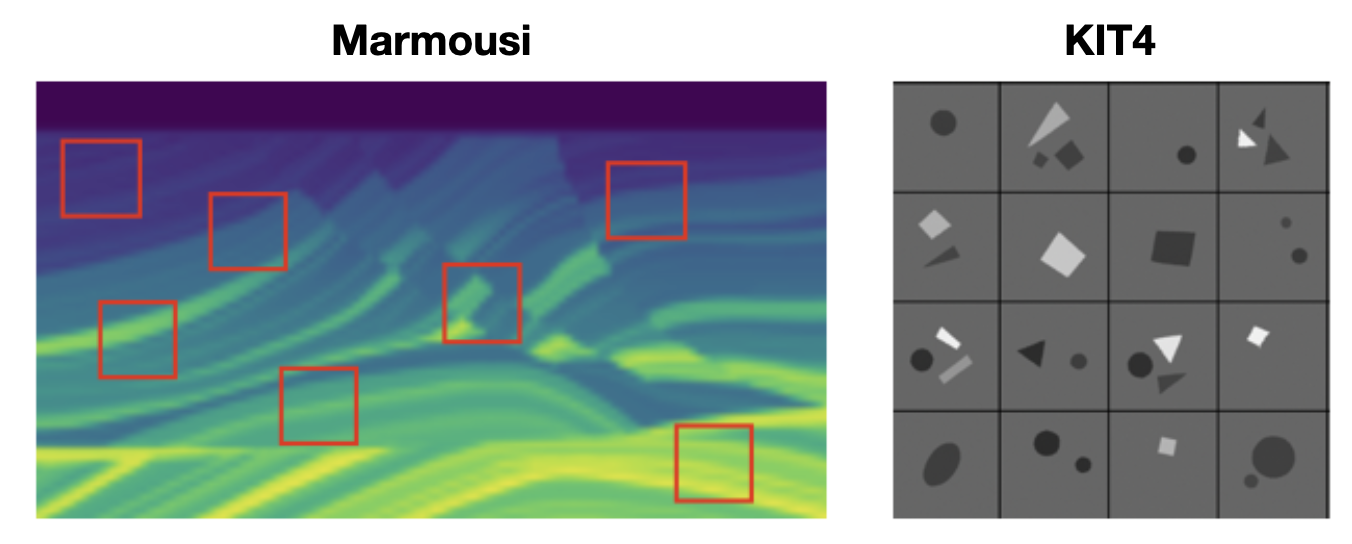}
    \caption{(a) The Marmousi dataset was collected by randomly extracting 128x128-sized sub-images from geological structure simulation images. (b) The KIT4 dataset was generated by randomly placing non-overlapping geometric shapes within a 128x128 area and assigning different speed field values to distinct regions.}
    \label{fig:data}
\end{figure}

In our numerical experiments, we consider two datasets: the Marmousi dataset \cite{martin2006marmousi2} and the KIT4 dataset \cite{hamilton2018deep} as illustrated in Fig.\ref{fig:data}. For the Marmousi dataset, we have constructed the dataset based on the simulated geological structure image to train the generative priors using diffusion models. Each pixel value in the $1000 \times 13600$ Marmousi image represents the wave propagation speed field within the structure. We randomly select $128 \times 128$ sub-images from this data, normalizing the values to the [0.01, 1] range. In this manner, we select $6000$ sub-images as the training dataset and $400$ sub-images as the testing dataset. The testing dataset is used not only for validating the score function during training but also for testing the subsequent travel-time tomography reconstruction problem. For the construction of the KIT4 dataset, we follow the approach described in \cite{cao2024dual}. This involves generating $1$ to $4$ regular geometric shapes and irregular polygon masks within the $[0,1] \times [0,1]$ region. These polygon masks are randomly rotated around the center at varying angles, ensuring that different polygon masks do not overlap. Each mask is then randomly assigned a speed field value within the [0.01, 1] range, while the background speed field is fixed at $0.5$. Following this, we similarly select $6000$ images as the training dataset and $400$ images as the testing dataset. 
  
\subsection{Receiver-Transmitter Geometry}
The reconstruction performance of inverse problems is strongly related to their inherent ill-posedness. For the considered travel-time tomography problem, only the first arrival time information is used, rather than the full-wave data employed in scattering tomography. Consequently, this inverse problem is expected to exhibit significant ill-posedness. Furthermore, the placement of receivers and transmitters on only a partial boundary makes the reconstruction even more challenging. From a practical application perspective, we consider three different geometries of receivers and transmitters across two simulated datasets.

\textbf{Horizontal Geometry.~} The first scenario involves oil exploration or geological surveys, where deep drilling is required to investigate subsurface cross-sectional structures. In our numerical example, we consider a $128 \times 128$ speed field $\x_0$ with horizontally opposed transmitters and receivers on the left and right sides. Specifically, $12$ transmitters are evenly spaced on the left side, and $24$ receivers are placed on the right side to collect first-arrival time information. Consequently, the observed data $\y$ in this scenario is a $24 \times 12$ matrix.

\textbf{Vertical Geometry.~} The second scenario involves reconstructing geological structures using the first-arrival times of seismic waves generated by natural earthquakes as they reach surface observation stations. In our numerical example, we consider a \(128 \times 128\) speed field $\x_0$ with corresponding transmitters and receivers on the top and bottom sides. Here, $12$ transmitters are equidistantly spaced on the bottom side, while $24$ receivers are placed on the top side to collect first-arrival time information. The observed data $\y$ in this scenario is a $24 \times 12$ matrix as well. 

\textbf{Surrounding  Geometry.~} The third scenario involves using transmitters and receivers placed around the surface of an object to detect its internal structure. In our numerical example, we consider a \(128 \times 128\) speed field $\x_0$ with the relatively dense placement of transmitters and receivers around all four sides. Specifically, 6 transmitters are evenly spaced on each of the four sides, and 24 receivers are evenly spaced on each of the four sides to collect first-arrival time information. Consequently, the observed data $\y$ in this scenario is a \(96 \times 24\) matrix.

\begin{figure}[!htp]
    \centering
    \includegraphics[width=0.9\textwidth]{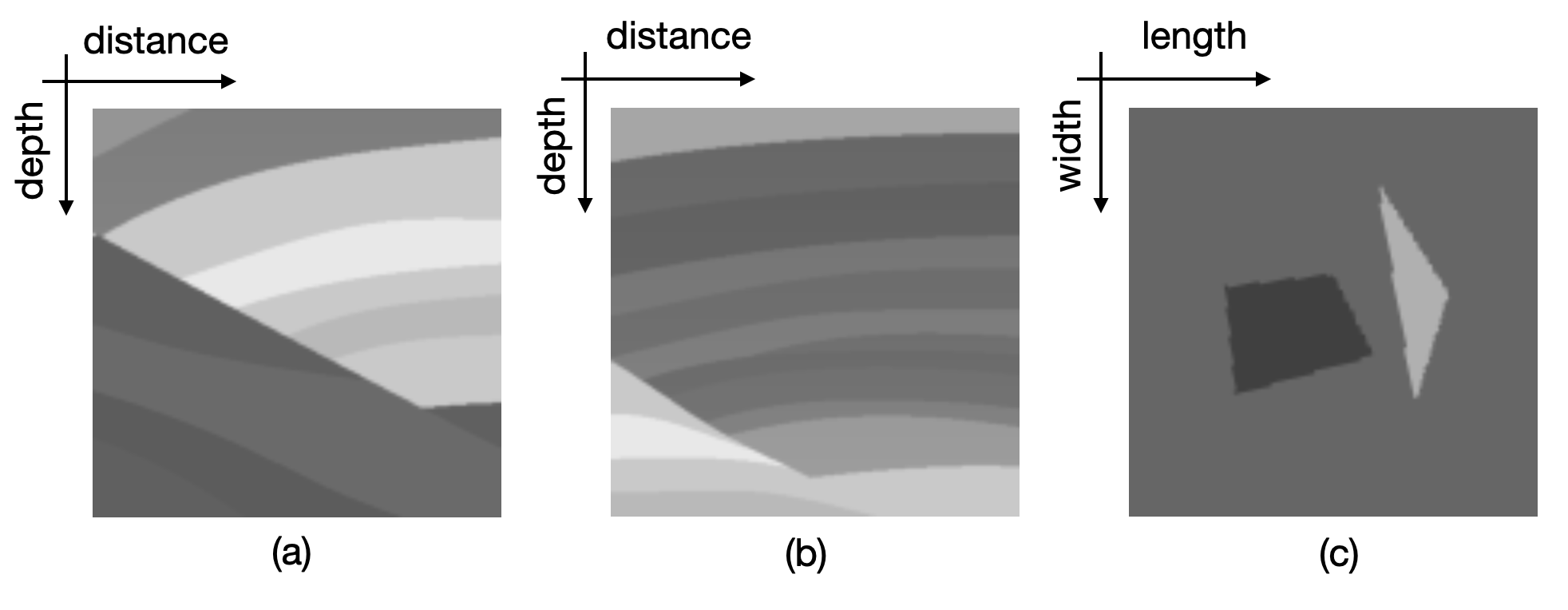}
    \caption{There are three geometries for placing the source-receiver pairs: (a)(b) For the Marmousi dataset, they are arranged in horizontal and vertical patterns; (c) For the KIT4 dataset, they are arranged in a circular pattern. As illustrated, the red pentagrams represent the locations of the sources, and the white triangles represent the locations of the receivers.}
    \label{mode}
\end{figure}

\subsection{Forward and Adjoint Modeling} 
The seismic wave is considered to be individually applied on different source points $\{(u^s_n, v^s_n)\}_{n=1}^{N}$, and the resulting travel time is measured at receiver point $\{(u^r_m, v^r_m)\}_{m=1}^{M}$ where the fast-marching method (FMM) \cite{sethian1999fast} is used for obtaining the discrete travel-time field $\y_n \in \Rd^{128\times 128}$ associated with the Eikonal equations. Specifically, the Eikonal equation is updated 
using the local solver with the upwind scheme:
\begin{equation}\label{local_solver_section_experiment}
   \sum\limits_{[a,b]\in \mathcal{N}_{[i,j]}} \left[\left(\frac{\y_n[i,j] - \y_n[a,b]}{h}\right)^{+}\right]^2 = \left(\frac{1}{\x_0[i,j]}\right)^2
\end{equation}
\begin{equation}\label{start_condition_section_discrete}
    \text{s.t. } \y_n\left[i^s_n, j^s_n\right] = 0
\end{equation}
where $\left[i^s_n, j^s_n\right] = \left[\lfloor \frac{u^s_n}{h} \rfloor,\lfloor\frac{v^s_n}{h}\rfloor\right]$ denotes the source index on $\x_0$, and $(x)^+$ is the positive part of $x$. Besides, FMM maintains a min-heap to update the wavefront propagation to ensure the point with the smallest travel time is processed first. Once the selected point with the minimal travel time in the wavefront is updated via \eqref{local_solver_section_experiment},  we pop out this point and append its neighbors into the heap until all grids have been updated. Furthermore, we consider the noised observation $\y^{\delta}_{m,n}$:
\begin{align}\label{forward_eikonal_numerical}
    \y_{m,n}^\delta :=  \y_{n}(u_m^r, v_m^r) + \delta \bm{\eta}_{m,n} \in \R,
\end{align}
where $\delta$ denotes the noise level of the simulated data, and $\bm{\eta}_{m,n}$ are independent noise samples from the standard Gaussian distribution. The complete data matrix $\y^{\delta} \in \R^{M\times N}$ is concatenated by $M \times N$ noised data $\y^{\delta}_{m,n}$, and the corresponding source locations $\{(u^s_n, v^s_n)\}_{n=1}^{N}$ are also known. 

To numerically calculate the adjoint-state variable $\boldsymbol{\Lambda}_n$ in \eqref{adjoint_eqn} \eqref{adjoint_eqn_cond}, we use the discretization $\boldsymbol{\lambda}_n^k$ for approximation 
\begin{equation}\label{adjoint_discrete_num}
  \sum\limits_{[a,b]\in \mathcal{N}_{[i,j]}}\left[ \left(\frac{\y_n[i,j] - \y_n[a,b]}{h}\right)^{+} \frac{\boldsymbol{\lambda}_n[i,j]}{h} - \left(\frac{\y_n[a,b] -\y_n[i,j]}{h}\right)^{+}\frac{\boldsymbol{\lambda}_n[a,b]}{h}\right] = 0,
\end{equation}
\begin{equation}\label{adjoint_discrete_cond_num}
     \text{s.t. } \left(\frac{\y_n[p,q] - \y_n[i^r_m,j^r_m]}{h}\right) \boldsymbol{\lambda}_n[i^r_m,j^r_m]  = \y^{\delta}_{m,n} - \y_n[i^r_m,j^r_m],
\end{equation}
where $\left[i^r_m, j^r_m\right] = \left[\lfloor \frac{u^r_m}{h} \rfloor,\lfloor\frac{v^r_m}{h}\rfloor\right]$ denotes the receiver index on $\x_0$, and $[p,q] \in \mathcal{N}_{[i,j]}$ is off the boundary. Note that the discrete adjoint equation \eqref{adjoint_discrete_num} \eqref{adjoint_discrete_cond_num} can be solved at each point following the reverse pop order provided by FMM \cite{deckelnick2011numerical}.

In this paper, the observed data matrix $\y^{\delta}  \in \R^{24 \times 12}$ in the horizontal and vertical geometries for the Marmousi dataset, while  $\y^{\delta}  \in \R^{96 \times 24}$ in the surrounding geometry for the KIT4 dataset. We also consider the simulated noise level $\delta = 2.5\%$ for all datasets.

\subsection{Training Subspace Score Models}
Here, we list the shared parameters for training the score models separately for all subspaces, as well as the divergence criterion to select the downsampling times during the forward diffusion process.

\textbf{Score Models.} We adopt almost the same parameter configurations as \cite{song2020score} for training the score functions. The model architecture is configured with a linear noise schedule with $\beta(t)$ values ranging from $0.1$ to $20$, and the time-embedding type is set to Fourier embeddings. Optimization is performed using the Adam optimizer with a learning rate of $2\times 10^{-4}$ for all layers, and gradient clipping is set at $1$ to prevent the gradient explosion. The training is performed with a batch size of $128$ across $1000$ epochs, and the optimal score model is selected with minimal validation loss. For the sampling process, we discretize the reverse SDE into $1000$ steps and use the Euler-Maruyama scheme for discretization. Note that we do not utilize the Predictor-Corrector (PC) algorithm, and instead perform an additional denoising step based on the Tweedie formula in the last step. Here, the $\frac{\epsilon}{T}$ is set to $10^{-3}$ for this denoising step. 

\begin{figure}[htbp]
    \centering
    \begin{minipage}{0.05\textwidth}
    \end{minipage}%
    \begin{minipage}{0.85\textwidth}
        \centering
        \includegraphics[width=\textwidth]{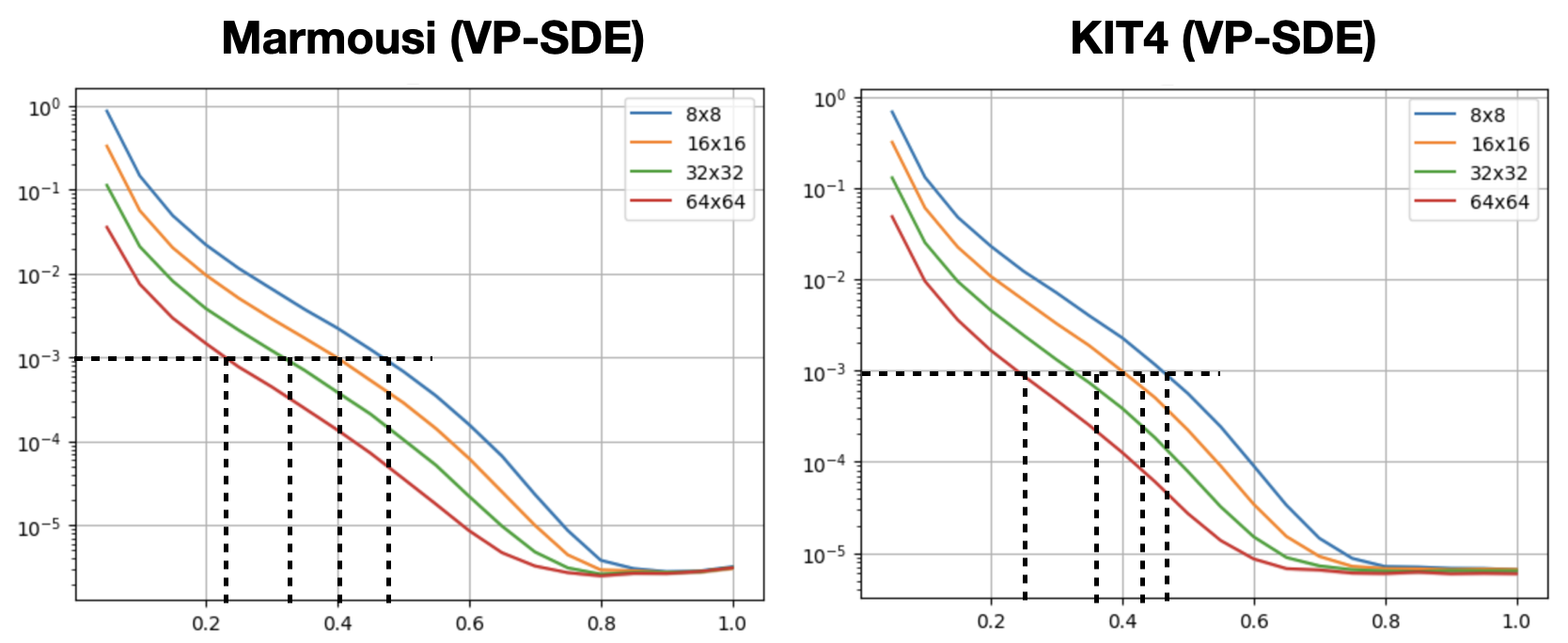} % 替换为你的图片文件
    \end{minipage}
    \caption{$D_F(\mathbf{U}_{k \mid k-1}; t)$ plots for two datasets are computed with respect to the trained full-space score model $\s_0(\x_t, t)$, where the horizontal axis represents the range of the sampling time, while the vertical axis indicates the values of the Orthogonal Fisher divergence. For a given divergence threshold $10^{-3}$, the optimal downsampling times $t_k$ for any subspace sequence are identified as the times when the corresponding divergences reach that threshold.}
    \label{fig:fisher}
\end{figure}

\textbf{Downsampling Times.} Starting with an original image size of $128\times 128$, we utilize the downsampling technique outlined earlier to reduce the dimensionality and train the score function $\s_k(\x_t^k, t)$ in each lower-dimensional space.
Besides, we adopt the method used in \cite{jing2022subspace} to determine the optimal downsampling times $t_k$, where the so-called orthogonal Fisher divergence is introduced 
\begin{equation}\label{eq:fisher}
    D_F(\mathbf{U}_{k \mid k-1}; t) = \frac{\Sigma^{k \mid k-1}_t}{d_{k-1} - d_k}\mathbb{E}_{\x_t}\left[\left\lVert \mathbf{P}^\perp_{k\mid k-1}\mathbf{U}_{k}^T\s_0(\x_t, t) + \frac{\x^{k \mid k-1}_t}{\Sigma^{k \mid k-1}_t} \right\rVert^2 \right]
\end{equation}
to quantify the divergence between the orthogonal components $\x_t^{k\mid k-1}$ and the Gaussian distribution with variance $\Sigma^{k \mid k-1}_t$. As shown in Fig.\ref{fig:fisher}, because the diffusion processes in intermediate dimensions are brief and the downsampled sequence with too many subspaces can degrade sampling quality, we limit our choice to the  $32 \times 32$ and $64 \times 64$ resolutions for subspace-based sampling, along with their corresponding downsampling times.

\section{Experiment Results}\label{6}
   
\subsection{Reconstruction Results}
We evaluate the reconstruction results of our proposed method on simulated test samples from the Marmousi and KIT4 datasets.  To demonstrate the numerical superiority of the proposed approach, we compare its performance with the following learning/non-learning methods. 
\begin{itemize}
    \item \textbf{L-BFGS.~} \cite{biondi2021object} leverages the linearization approach incorporated with the L-BFGS optimization algorithm for gradient update. In our experiments, we choose the maximum iteration to be $30$ to guarantee the convergence and apply the Gaussian smoothing to the computed gradient to avoid singularities during the optimization process.

    \item \textbf{DIP.~} The DIP-based method employs the Deep Image Prior \cite{ulyanov2018deep} as an implicit regularization for reconstructing the velocity field. In a similar vein as \cite{liu2023deepeit} for solving inverse problems, we use an untrained U-net to re-parameterize the velocity field with the fixed noise as network input. Early stopping is utilized as a regularization technique in this method.    

    \item \textbf{RED-Diff.~} The RED-Diff \cite{mardani2023variational} method samples from a posterior distribution based on the well-trained diffusion prior and formulates the posterior sampling as a regularized stochastic optimization. The regularization is imposed by the denoising diffusion process with denoisers at different timesteps, while the data discrepancy term is updated by the adjoint-state \cite{leung2006adjoint} method.   

    \item \textbf{Supervised Post-Processing.~} The learning-based post-processing method \cite{araya2018deep} is widely used for solving inverse problems in an end-to-end fashion. In our experiments, the L-BFGS algorithm produces a rough approximation for the velocity reconstruction as the network input. The network parameters are trained to match the ground truth as the output. Unlike the aforementioned comparison methods, this supervised method relies on amounts of paired noisy measurements and samples for supervised learning, which presents challenges when a large training dataset is not available.

\end{itemize}
Note that the supervised method typically outperforms zero-shot learning approaches, as the availability of large paired datasets enhances the generalization ability of the supervised methods. In contrast, the other three approaches and our proposed methods do not leverage this prior knowledge.

\begin{figure}[!htp]
    \centering
    \includegraphics[width=0.9\textwidth]{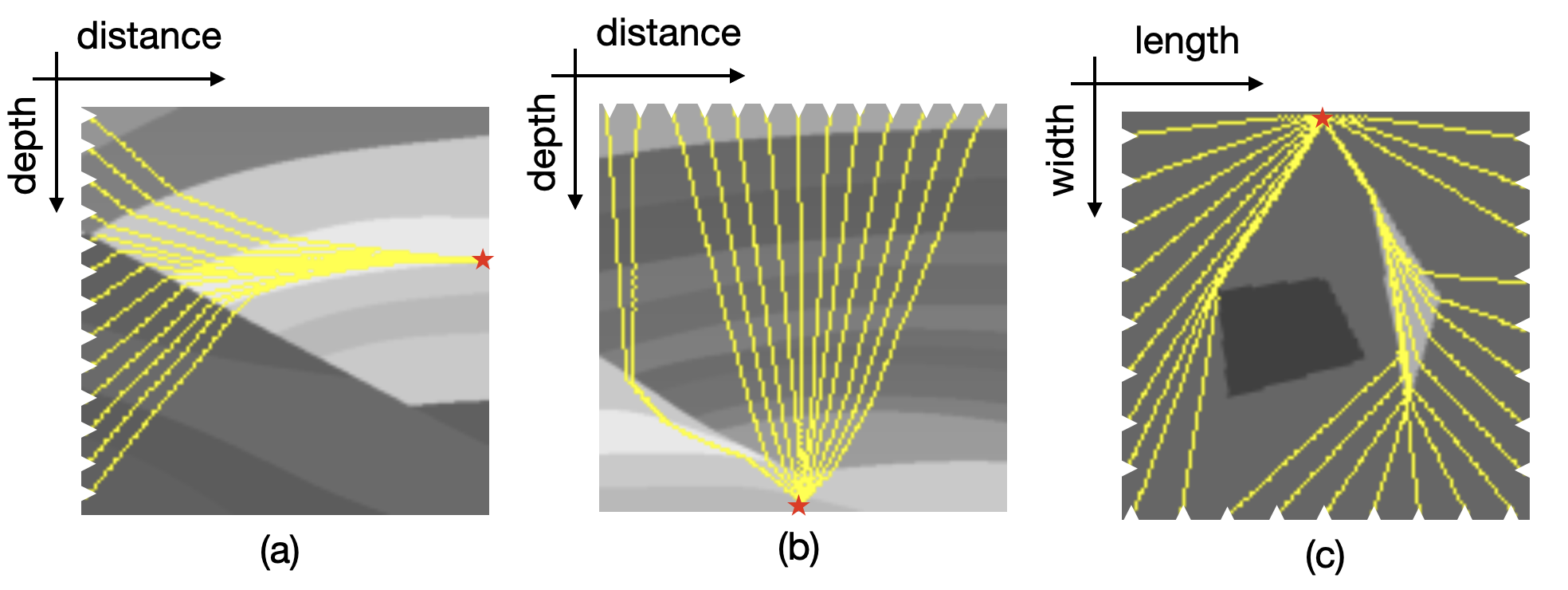}
    \caption{The ray paths travel from each transmitter to the receivers, following the route that minimizes the travel time. (a) In the horizontal geometry, the ray paths are prone to traverse the high-speed fault region to reach the receivers. (b) In the vertical geometry, the ray paths encounter different faults with different speed values. (c) In the surrounding pattern, the rays tend to bypass the subregions with lower speed values and instead pass through higher-speed areas. }
    \label{trace1}
\end{figure}

\begin{figure}[!htp]
  \centering
  \scalebox{0.9}{
  \begin{tabular}{cccccccccccccc}
    33&\scriptsize{\makecell{L-BFGS}}&\scriptsize{DIP }&
    \scriptsize{RED-Diff}&\scriptsize{\makecell{Supervised \\ Post-Processing}}&\scriptsize{\makecell{PDE-based \\ Full-space DPS}}&\scriptsize{\makecell{PDE-based \\ Subspace DPS}}&\scriptsize{\makecell{Ground \\ Truth}}
    \\ 
    \put(-5,10){\rotatebox{90}{\scriptsize{Sample 1}}}&
    \includegraphics[width=.12\linewidth,height=.12\linewidth]{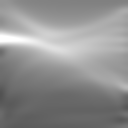}&
    \includegraphics[width=.12\linewidth,height=.12\linewidth]{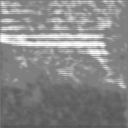}&
    \includegraphics[width=.12\linewidth,height=.12\linewidth]{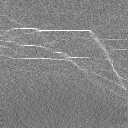}&
    \includegraphics[width=.12\linewidth,height=.12\linewidth]{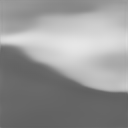}&
    \includegraphics[width=.12\linewidth,height=.12\linewidth]{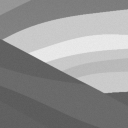}&
    \includegraphics[width=.12\linewidth,height=.12\linewidth]{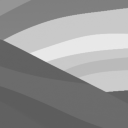}&
    \includegraphics[width=.12\linewidth,height=.12\linewidth]{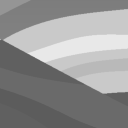}
    \\ 
    \put(-5,10){\rotatebox{90}{\scriptsize{Sample 2}}}&
    \includegraphics[width=.12\linewidth,height=.12\linewidth]{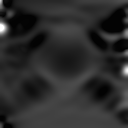}&
    \includegraphics[width=.12\linewidth,height=.12\linewidth]{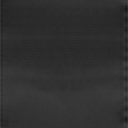}&
    \includegraphics[width=.12\linewidth,height=.12\linewidth]{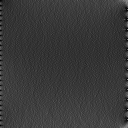}&
    \includegraphics[width=.12\linewidth,height=.12\linewidth]{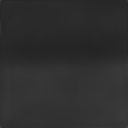}&
    \includegraphics[width=.12\linewidth,height=.12\linewidth]{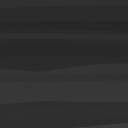}&
    \includegraphics[width=.12\linewidth,height=.12\linewidth]{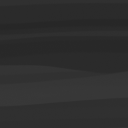}&
    \includegraphics[width=.12\linewidth,height=.12\linewidth]{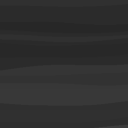}
    \\ 
    \put(-5,10){\rotatebox{90}{\scriptsize{Sample 3}}}&
    \includegraphics[width=.12\linewidth,height=.12\linewidth]{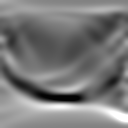}&
    \includegraphics[width=.12\linewidth,height=.12\linewidth]{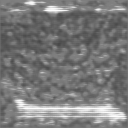}&
    \includegraphics[width=.12\linewidth,height=.12\linewidth]{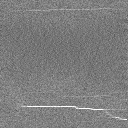}&
    \includegraphics[width=.12\linewidth,height=.12\linewidth]{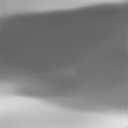}&
    \includegraphics[width=.12\linewidth,height=.12\linewidth]{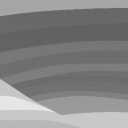}&
    \includegraphics[width=.12\linewidth,height=.12\linewidth]{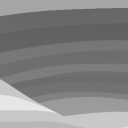}&
    \includegraphics[width=.12\linewidth,height=.12\linewidth]{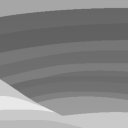}
    \\
     \put(-5,10){\rotatebox{90}{\scriptsize{Sample 4}}}&
   \includegraphics[width=.12\linewidth,height=.12\linewidth]{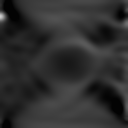}&
    \includegraphics[width=.12\linewidth,height=.12\linewidth]{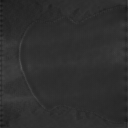}&
    \includegraphics[width=.12\linewidth,height=.12\linewidth]{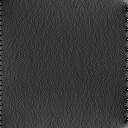}&
    \includegraphics[width=.12\linewidth,height=.12\linewidth]{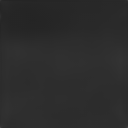}&
    \includegraphics[width=.12\linewidth,height=.12\linewidth]{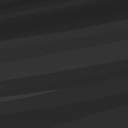}&
    \includegraphics[width=.12\linewidth,height=.12\linewidth]{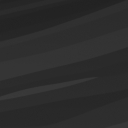}&
    \includegraphics[width=.12\linewidth,height=.12\linewidth]{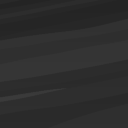}

    \\
     \put(-5,10){\rotatebox{90}{\scriptsize{Sample 5}}}&
    \includegraphics[width=.12\linewidth,height=.12\linewidth]{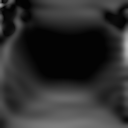}&
    \includegraphics[width=.12\linewidth,height=.12\linewidth]{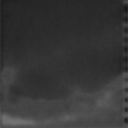}&
    \includegraphics[width=.12\linewidth,height=.12\linewidth]{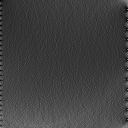}&
    \includegraphics[width=.12\linewidth,height=.12\linewidth]{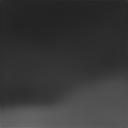}&
    \includegraphics[width=.12\linewidth,height=.12\linewidth]{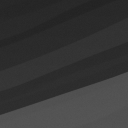}&
    \includegraphics[width=.12\linewidth,height=.12\linewidth]{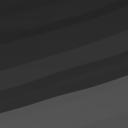}&
    \includegraphics[width=.12\linewidth,height=.12\linewidth]{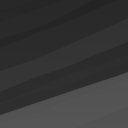}
    \\
     \put(-5,10){\rotatebox{90}{\scriptsize{Sample 6}}}&
   \includegraphics[width=.12\linewidth,height=.12\linewidth]{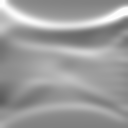}&
    \includegraphics[width=.12\linewidth,height=.12\linewidth]{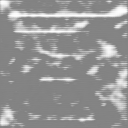}&
    \includegraphics[width=.12\linewidth,height=.12\linewidth]{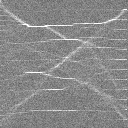}&
    \includegraphics[width=.12\linewidth,height=.12\linewidth]{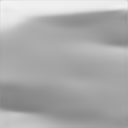}&
    \includegraphics[width=.12\linewidth,height=.12\linewidth]{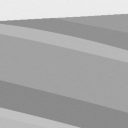}&
    \includegraphics[width=.12\linewidth,height=.12\linewidth]{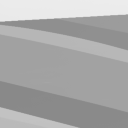}&
    \includegraphics[width=.12\linewidth,height=.12\linewidth]{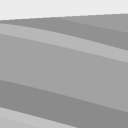}
  \end{tabular}}
  \caption{The reconstruction results for each Marmousi test sample are compared across comparison methods and ours, using measurements from horizontally placed sources and receivers.}
  \label{Marmousi_1}
\end{figure}

Fig.\ref{Marmousi_1} shows a visual comparison of reconstruction results for horizontally placed transmitter-receiver pairs on the Marmousi dataset. In this configuration, the ray paths in Fig.\ref{trace1}(a) originating from the source almost traverse through fault zones with higher velocity values in the nearly horizontal direction. Because the fault velocity values along this path are consistent, the inverse problem associated with this measurement approach exhibits relatively weak ill-posedness. It can be seen that the proposed PDE-based full-space and subspace DPS methods significantly outperform the four comparison methods, almost both achieving a perfect recovery of the underlying parameter distribution. The L-BFGS method resorts to the linearization technique and the Gaussian smoothing for minimizing the PDE-based measurement misfit. It yields quite blurred reconstructions due to the smoothed gradients and the lack of prior knowledge. In the DIP method, the required step size and the stopping time are manually tunable parameters. Even when we select the optimal parameters for optimization, the reconstruction results still contain amounts of misleading artifacts. Another optimization method RED-diff employs score models as prior knowledge, however, fails to recover the velocity field. Compared to the supervised post-processing method, the PDE-based full-space/subspace DPS approaches yield clearer reconstructed images and provide more accurate reconstruction of interfaces between different faults, even though our methods do not rely on amounts of paired data for training.    

\begin{figure}[!htp]
  \centering
  \scalebox{0.9}{
  \begin{tabular}{cccccccccccccc}
    %\hline
    %Items & Advantages & Disadvantages
    &\scriptsize{\makecell{L-BFGS}}&\scriptsize{DIP }&
    \scriptsize{RED-Diff}&\scriptsize{\makecell{Supervised \\ Post-Processing}}&\scriptsize{\makecell{PDE-based \\ Full-space DPS}}&\scriptsize{\makecell{PDE-based \\ Subspace DPS}}&\scriptsize{\makecell{Ground \\ Truth}}
    \\ 
    \put(-5,10){\rotatebox{90}{\scriptsize{Sample 1}}}&
    \includegraphics[width=.12\linewidth,height=.12\linewidth]{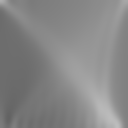}&
    \includegraphics[width=.12\linewidth,height=.12\linewidth]{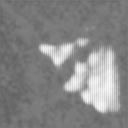}&
    \includegraphics[width=.12\linewidth,height=.12\linewidth]{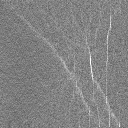}&
    \includegraphics[width=.12\linewidth,height=.12\linewidth]{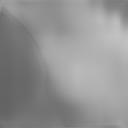}&
    \includegraphics[width=.12\linewidth,height=.12\linewidth]{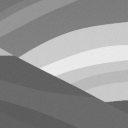}&
    \includegraphics[width=.12\linewidth,height=.12\linewidth]{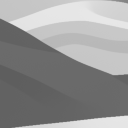}&
    \includegraphics[width=.12\linewidth,height=.12\linewidth]{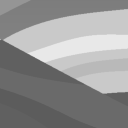}
    \\ 
    \put(-5,10){\rotatebox{90}{\scriptsize{Sample 2}}}&
    \includegraphics[width=.12\linewidth,height=.12\linewidth]{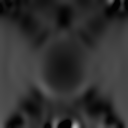}&
    \includegraphics[width=.12\linewidth,height=.12\linewidth]{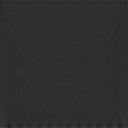}&
    \includegraphics[width=.12\linewidth,height=.12\linewidth]{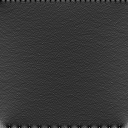}&
    \includegraphics[width=.12\linewidth,height=.12\linewidth]{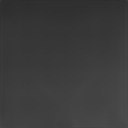}&
    \includegraphics[width=.12\linewidth,height=.12\linewidth]{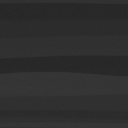}&
    \includegraphics[width=.12\linewidth,height=.12\linewidth]{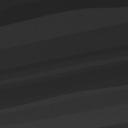}&
    \includegraphics[width=.12\linewidth,height=.12\linewidth]{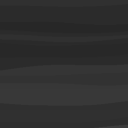}
    \\ 
    \put(-5,10){\rotatebox{90}{\scriptsize{Sample 3}}}&
    \includegraphics[width=.12\linewidth,height=.12\linewidth]{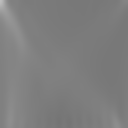}&
    \includegraphics[width=.12\linewidth,height=.12\linewidth]{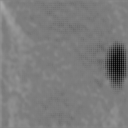}&
    \includegraphics[width=.12\linewidth,height=.12\linewidth]{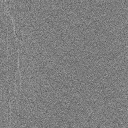}&
    \includegraphics[width=.12\linewidth,height=.12\linewidth]{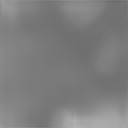}&
    \includegraphics[width=.12\linewidth,height=.12\linewidth]{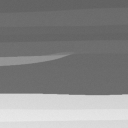}&
    \includegraphics[width=.12\linewidth,height=.12\linewidth]{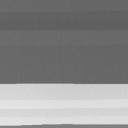}&
    \includegraphics[width=.12\linewidth,height=.12\linewidth]{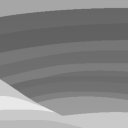}
    \\
     \put(-5,10){\rotatebox{90}{\scriptsize{Sample 4}}}&
     \includegraphics[width=.12\linewidth,height=.12\linewidth]{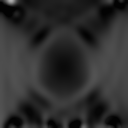}&
    \includegraphics[width=.12\linewidth,height=.12\linewidth]{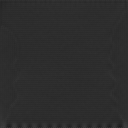}&
    \includegraphics[width=.12\linewidth,height=.12\linewidth]{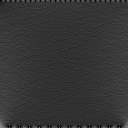}&
    \includegraphics[width=.12\linewidth,height=.12\linewidth]{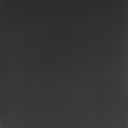}&
    \includegraphics[width=.12\linewidth,height=.12\linewidth]{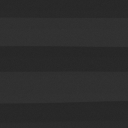}&
    \includegraphics[width=.12\linewidth,height=.12\linewidth]{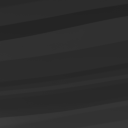}&
    \includegraphics[width=.12\linewidth,height=.12\linewidth]{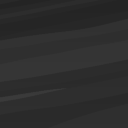}

    \\
     \put(-5,10){\rotatebox{90}{\scriptsize{Sample 5}}}&
    \includegraphics[width=.12\linewidth,height=.12\linewidth]{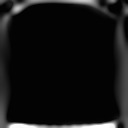}&
    \includegraphics[width=.12\linewidth,height=.12\linewidth]{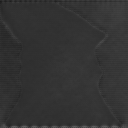}&
    \includegraphics[width=.12\linewidth,height=.12\linewidth]{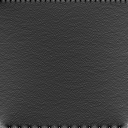}&
    \includegraphics[width=.12\linewidth,height=.12\linewidth]{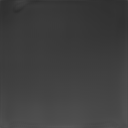}&
    \includegraphics[width=.12\linewidth,height=.12\linewidth]{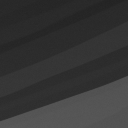}&
    \includegraphics[width=.12\linewidth,height=.12\linewidth]{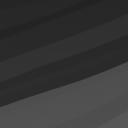}&
    \includegraphics[width=.12\linewidth,height=.12\linewidth]{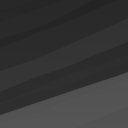}
    \\
     \put(-5,10){\rotatebox{90}{\scriptsize{Sample 6}}}&
    \includegraphics[width=.12\linewidth,height=.12\linewidth]{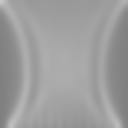}&
    \includegraphics[width=.12\linewidth,height=.12\linewidth]{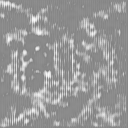}&
    \includegraphics[width=.12\linewidth,height=.12\linewidth]{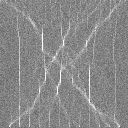}&
    \includegraphics[width=.12\linewidth,height=.12\linewidth]{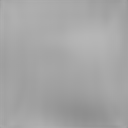}&
    \includegraphics[width=.12\linewidth,height=.12\linewidth]{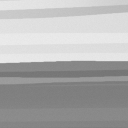}&
    \includegraphics[width=.12\linewidth,height=.12\linewidth]{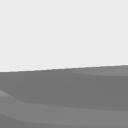}&
    \includegraphics[width=.12\linewidth,height=.12\linewidth]{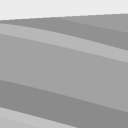}
  \end{tabular}}
  \caption{The reconstruction results for each Marmousi test sample are compared across comparison methods and ours, using measurements from vertically placed sources and receivers.}
  \label{Marmousi_2}
\end{figure}

Another visual comparison of numerical results is depicted in Fig.\ref{Marmousi_2} for vertically placed transmitter-receiver pairs. The waves propagate through various faults and reach the receiver located near the surface in Fig.\ref{trace1}(b), where each ray path will encounter faults with significantly different velocity values. Meanwhile, small differences in measurements at adjacent receivers make the inverse problem more prone to strong ill-posedness. As shown in Fig.\ref{Marmousi_2}, none of the other comparison methods can produce reasonable reconstruction results that align with the underlying spatial characteristics, while our proposed method yields more satisfactory estimates for the velocity field with sharp fault boundaries. Notice that, the reconstructed samples still exhibit some differences from the ground truths due to the ill-posedness of the problem, whereas the boundary measurement misfit has converged to a small value.   

\begin{figure}[!htp]
  \centering
  \scalebox{0.9}{
  \begin{tabular}{cccccccccccccc}
    %\hline
    %Items & Advantages & Disadvantages
    &\scriptsize{\makecell{L-BFGS}}&\scriptsize{DIP }&
    \scriptsize{RED-Diff}&\scriptsize{\makecell{Supervised \\ Post-Processing}}&\scriptsize{\makecell{PDE-based \\ Full-space DPS}}&\scriptsize{\makecell{PDE-based \\ Subspace DPS}}&\scriptsize{\makecell{Ground \\ Truth}}
    \\ 
    \put(-5,10){\rotatebox{90}{\scriptsize{Sample 1}}}&
    \includegraphics[width=.12\linewidth,height=.12\linewidth]{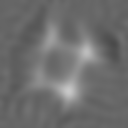}&
    \includegraphics[width=.12\linewidth,height=.12\linewidth]{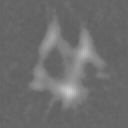}&
    \includegraphics[width=.12\linewidth,height=.12\linewidth]{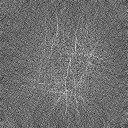}&
    \includegraphics[width=.12\linewidth,height=.12\linewidth]{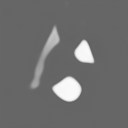}&
    \includegraphics[width=.12\linewidth,height=.12\linewidth]{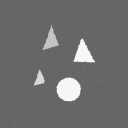}&
    \includegraphics[width=.12\linewidth,height=.12\linewidth]{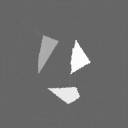}&
    \includegraphics[width=.12\linewidth,height=.12\linewidth]{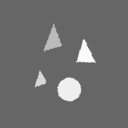}
    \\ 
    \put(-5,10){\rotatebox{90}{\scriptsize{Sample 2}}}&
    \includegraphics[width=.12\linewidth,height=.12\linewidth]{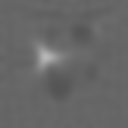}&
    \includegraphics[width=.12\linewidth,height=.12\linewidth]{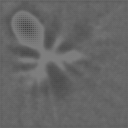}&
    \includegraphics[width=.12\linewidth,height=.12\linewidth]{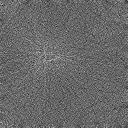}&
    \includegraphics[width=.12\linewidth,height=.12\linewidth]{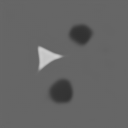}&
    \includegraphics[width=.12\linewidth,height=.12\linewidth]{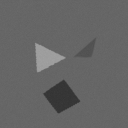}&
    \includegraphics[width=.12\linewidth,height=.12\linewidth]{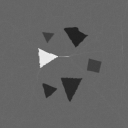}&
    \includegraphics[width=.12\linewidth,height=.12\linewidth]{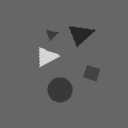}
    \\
     \put(-5,10){\rotatebox{90}{\scriptsize{Sample 3}}}&
   \includegraphics[width=.12\linewidth,height=.12\linewidth]{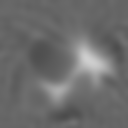}&
    \includegraphics[width=.12\linewidth,height=.12\linewidth]{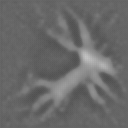}&
    \includegraphics[width=.12\linewidth,height=.12\linewidth]{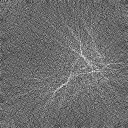}&
    \includegraphics[width=.12\linewidth,height=.12\linewidth]{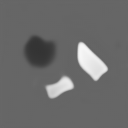}&
    \includegraphics[width=.12\linewidth,height=.12\linewidth]{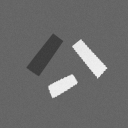}&
    \includegraphics[width=.12\linewidth,height=.12\linewidth]{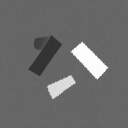}&
    \includegraphics[width=.12\linewidth,height=.12\linewidth]{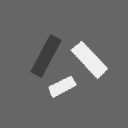}
    \\ 
    \put(-5,10){\rotatebox{90}{\scriptsize{Sample 4}}}&
    \includegraphics[width=.12\linewidth,height=.12\linewidth]{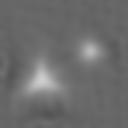}&
    \includegraphics[width=.12\linewidth,height=.12\linewidth]{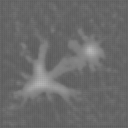}&
    \includegraphics[width=.12\linewidth,height=.12\linewidth]{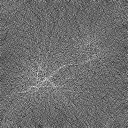}&
    \includegraphics[width=.12\linewidth,height=.12\linewidth]{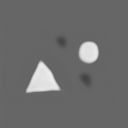}&
    \includegraphics[width=.12\linewidth,height=.12\linewidth]{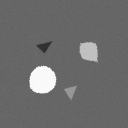}&
    \includegraphics[width=.12\linewidth,height=.12\linewidth]{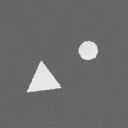}&
    \includegraphics[width=.12\linewidth,height=.12\linewidth]{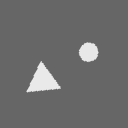}

    \\
     \put(-5,10){\rotatebox{90}{\scriptsize{Sample 5}}}&
    \includegraphics[width=.12\linewidth,height=.12\linewidth]{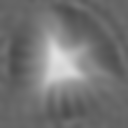}&
    \includegraphics[width=.12\linewidth,height=.12\linewidth]{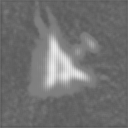}&
    \includegraphics[width=.12\linewidth,height=.12\linewidth]{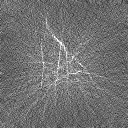}&
    \includegraphics[width=.12\linewidth,height=.12\linewidth]{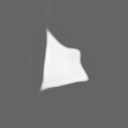}&
    \includegraphics[width=.12\linewidth,height=.12\linewidth]{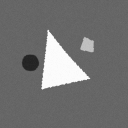}&
    \includegraphics[width=.12\linewidth,height=.12\linewidth]{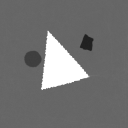}&
    \includegraphics[width=.12\linewidth,height=.12\linewidth]{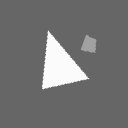}
    \\
     \put(-5,10){\rotatebox{90}{\scriptsize{Sample 6}}}&
    \includegraphics[width=.12\linewidth,height=.12\linewidth]{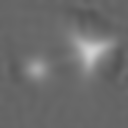}&
    \includegraphics[width=.12\linewidth,height=.12\linewidth]{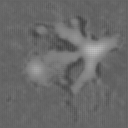}&
    \includegraphics[width=.12\linewidth,height=.12\linewidth]{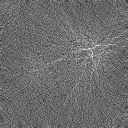}&
    \includegraphics[width=.12\linewidth,height=.12\linewidth]{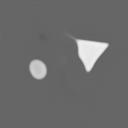}&
    \includegraphics[width=.12\linewidth,height=.12\linewidth]{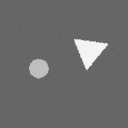}&
    \includegraphics[width=.12\linewidth,height=.12\linewidth]{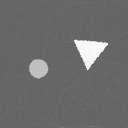}&
    \includegraphics[width=.12\linewidth,height=.12\linewidth]{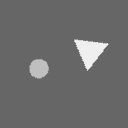}
  \end{tabular}}
  \caption{The reconstruction results for each KIT4 test sample are compared across comparison methods and ours, using measurements from sources and receivers placed around the object.}
  \label{KIT4_1}
\end{figure}

For the KIT4 dataset, the transmitters and receivers are placed around the considered region. Since the ray path of the propagating wave in Fig.\ref{trace1}(c) tends to bend around low-velocity subregions, the first arrival time associated with each ray path lacks the velocity information from these regions, making it difficult to accurately reconstruct these areas. Fig.\ref{KIT4_1} plots the reconstruction results in each row for different methods. We observe that the first three methods are capable of reconstructing the rough contours of the internal subregions, However, they often generate unexpected artifacts and result in incorrect velocity estimations. The supervised method heavily relies on rough approximations as input, while our proposed methods sample the reconstruction starting from the random noise and produce better shape reconstruction and velocity estimation.

\subsection{Quantitative Results}
For quantitative comparison, we select two widely used metrics: \textbf{RMSE} (Root Mean Square Error) and \textbf{SSIM} (Structural Similarity Index) for quality assessment. A higher SSIM and a lower RMSE signify a closer match between the reconstruction and the ground truth. 

In Table.\ref{quantitative_comp}, we calculate the average RMSE and SSIM metrics across two test datasets, each containing 400 test samples. Compared to the three unsupervised methods, our proposed method achieves the best performance for three different configurations. Note that in the vertical configuration for the Marmousi dataset and the surrounding configuration for the KIT4 dataset, the supervised post-processing method outperforms our approach in terms of quantitative metrics. This also indicates that a large amount of measurement and sample data pairs can indeed provide strong prior knowledge for learning-based methods in solving inverse problems.     

\begin{table}[t]
\centering
\resizebox{\textwidth}{!}{% <------ Don't forget this %
\begin{tabular}{cccccccc}
\toprule
{} & \multicolumn{2}{c}{\textbf{Marmousi (horizontal)}} & \multicolumn{2}{c}{\textbf{Marmousi (vertical)}} & \multicolumn{2}{c}{\textbf{KIT4 (surround)}} \\
\cmidrule(lr){2-3}
\cmidrule(lr){4-5}
\cmidrule(lr){6-7}

{\thead{\textbf{Method}}} & {RMSE $\downarrow$} & {SSIM $\uparrow$} & {RMSE $\downarrow$} & {SSIM $\uparrow$} & {RMSE $\downarrow$} & {SSIM $\uparrow$} \\  

\midrule
\thead{L-BFGS} & 0.128 & 0.457 & 0.166 & 0.359 & 0.054 & 0.713 & \\
\thead{DIP} & 0.147 & 0.351 &  0.166 & 0.264 & 0.066 & 0.342 & \\
\thead{RED-Diff} & 0.173 & 0.193 & 0.178 & 0.168 &  0.105 & 0.085 & \\
\thead{Supervised Post-Processing} & 0.039 & 0.655 & \underline{0.104} & \underline{0.491} & \textbf{0.031} & \textbf{0.925} & \\ 
\cmidrule(l){1-7}  
\thead{PDE-based Full-space DPS} & \textbf{0.014} & \textbf{0.768}  & \textbf{0.089} & \textbf{0.525} &  \underline{0.050} & \underline{0.864}&\\
\thead{PDE-based Subspace DPS} & \underline{0.024} & \underline{0.679} & 0.113 & 0.462 & 0.062 & 0.837 & \\  
     
\bottomrule
\end{tabular}
}
\caption{The quantitative results (RMSE, SSIM) are evaluated on the Marmousi and KIT4 test datasets for different  methods. \textbf{Bold}: best, \underline{underline}: second best.
}
\label{quantitative_comp}
\end{table}

\subsection{Computational Efficiency}
One key motivation behind our proposed method is to improve the computational efficiency when applying the iterative sampling-based DPS approach to solving PDE-based inverse problems. The subspace-based technique not only reduces the dimensionality of score function evaluations but also decreases the number of discretization grids required by numerical PDE solvers.   

To provide a comprehensive comparison, we compare the sampling times associated with two downsampling strategies ($32 \rightarrow 64 \rightarrow 128$, $64 \rightarrow 128$) against those in the full-space sampling approach (See Table \ref{computational_effi}). There is a trade-off between reducing sampling time and improving the quality of the generated samples, that is, increasing the downsampling numbers will deteriorate the quality of the generated samples. Note that, the improvement in sampling acceleration comes from two important factors: the increased speed of score functions evaluated on the GPU and numerical PDE solvers executed on the CPU.
  
\begin{table}[t]
\centering
\resizebox{\textwidth}{!}{% <------ Don't forget this %
\begin{tabular}{lllllll}
\toprule
{} & \multicolumn{2}{c}{\textbf{KIT4 (surround)}} & \multicolumn{2}{c}{\textbf{Marmousi (vertical)}} & \multicolumn{2}{c}{\textbf{Marmousi (horizontal)}} \\
\cmidrule(lr){2-3}
\cmidrule(lr){4-5} 
\cmidrule(lr){6-7}

{\thead{\textbf{Method}}} & {$32 \rightarrow 64 \rightarrow 128$} & {$64 \rightarrow 128$} & {$32 \rightarrow 64 \rightarrow 128$} & {$64 \rightarrow 128$} & {$32 \rightarrow 64 \rightarrow 128$} & {$64 \rightarrow 128$} \\

\midrule
\thead{PDE-based Subspace DPS \\ Average Time (sec)} & \thead{351.06} & \thead{583.01} & \thead{193.39} & \thead{301.11} & \thead{190.12} & \thead{293.58} \\
\thead{PDE-based Full-space DPS \\ Average Time (sec)} & \thead{(1060.75)} & \thead{(1060.75)} & \thead{(551.73)} & \thead{(551.73)} & \thead{(548.46)} & \thead{(548.46)} \\

\midrule
\thead{Average Time Ratio} & \thead{0.331} & \thead{0.550} & \thead{0.350} & \thead{0.546} & \thead{0.346} & \thead{0.535} \\
\bottomrule
\end{tabular}
}
\caption{The computational times of the PDE-based full-space and subspace DPS methods are compared on the Marmousi and KIT4 test datasets, respectively.     
}
\label{computational_effi}   
\vspace{-0.6cm}
\end{table}

\textbf{Remarks.} All nonlinear PDE solvers for different transmitter-receiver pairs are parallelly executed on the Intel-Xeon-Platinum-8358P CPU, while the training and inference processes of the conditional score-based models are performed on the NVIDIA-A100 GPU.

\subsection{Effect of the step size}
\label{step_size_alpha}
\begin{figure}[!htp]
  \centering
  \scalebox{0.9}{
  \begin{tabular}{cccccccccc}
    %\hline
    %Items & Advantages & Disadvantages
    &\scriptsize{\makecell{Full-space DPS \\ (fixed step size)}}&\scriptsize{\makecell{Full-space DPS  \\(varying step size)}}& \scriptsize{\makecell{Subspace DPS \\ (fixed step size)}}&\scriptsize{\makecell{Subspace DPS \\(varying step size)}}&\scriptsize{\makecell{Ground \\ Truth}}
    \\ 
    \put(-5,10){\rotatebox{90}{\scriptsize{Sample 1}}}&
    \includegraphics[width=.12\linewidth,height=.12\linewidth]{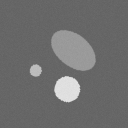}&
    \includegraphics[width=.12\linewidth,height=.12\linewidth]{TravelTime_DPS/KIT4_s/DPS_5.png}&
    \includegraphics[width=.12\linewidth,height=.12\linewidth]{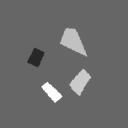}&
    \includegraphics[width=.12\linewidth,height=.12\linewidth]{TravelTime_DPS/KIT4_s/Subspace_5.png}&
    \includegraphics[width=.12\linewidth,height=.12\linewidth]{TravelTime_DPS/KIT4_s/GT_5.png}
    \\ 
    \put(-5,10){\rotatebox{90}{\scriptsize{Sample 1}}}&
    \includegraphics[width=.12\linewidth,height=.12\linewidth]{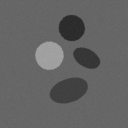}&
    \includegraphics[width=.12\linewidth,height=.12\linewidth]{TravelTime_DPS/KIT4_s/DPS_9.png}&
    \includegraphics[width=.12\linewidth,height=.12\linewidth]{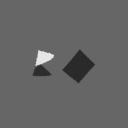}&
    \includegraphics[width=.12\linewidth,height=.12\linewidth]{TravelTime_DPS/KIT4_s/Subspace_9.png}&
    \includegraphics[width=.12\linewidth,height=.12\linewidth]{TravelTime_DPS/KIT4_s/GT_9.png}
    \\ 
    \put(-5,10){\rotatebox{90}{\scriptsize{Sample 1}}}&
    \includegraphics[width=.12\linewidth,height=.12\linewidth]{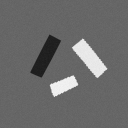}&
    \includegraphics[width=.12\linewidth,height=.12\linewidth]{TravelTime_DPS/KIT4_s/DPS_11.png}&
    \includegraphics[width=.12\linewidth,height=.12\linewidth]{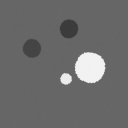}&
    \includegraphics[width=.12\linewidth,height=.12\linewidth]{TravelTime_DPS/KIT4_s/Subspace_11.png}&
    \includegraphics[width=.12\linewidth,height=.12\linewidth]{TravelTime_DPS/KIT4_s/GT_11.png}
    \\    
    \put(-5,10){\rotatebox{90}{\scriptsize{Sample 1}}}&
    \includegraphics[width=.12\linewidth,height=.12\linewidth]{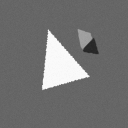}&
    \includegraphics[width=.12\linewidth,height=.12\linewidth]{TravelTime_DPS/KIT4_s/DPS_29.png}&
    \includegraphics[width=.12\linewidth,height=.12\linewidth]{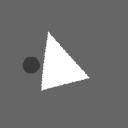}&   
    \includegraphics[width=.12\linewidth,height=.12\linewidth]{TravelTime_DPS/KIT4_s/Subspace_29.png}&
    \includegraphics[width=.12\linewidth,height=.12\linewidth]{TravelTime_DPS/KIT4_s/GT_29.png}
    \\
  \end{tabular}}
  \caption{The reconstruction results of the full-space/subspace DPS methods are compared when different step size schedules are applied.}
  \label{step_size}
\end{figure}

As noted by \cite{chung2023diffusion}, the step size $\rho'$ plays a crucial role in weighting the data consistency of the inverse problem. The authors suggest using $\rho':= \rho\left\|\y - \Ac(\hat\x_0(\x_t))\right\|$ as a step size schedule to effectively balance the measurement misfit with the score prior. For this varying step size schedule, the adaptive step size $\rho'$ is adjusted based on the data discrepancy $\left\|\y - \Ac(\hat\x_0(\x_t))\right\|$. Therefore, a larger step size ensures that the early stage of the reverse sampling process is more inclined to guarantee data consistency, while the later stage becomes more focused on gradual denoising to generate samples. As shown in Fig.\ref{step_size}, compared to the fixed step size schedule, using the varying step size schedule in PDE-based full-space/subspace DPS can improve the accuracy of inverse problem reconstruction.

 \section{Conclusion}\label{7}
In this paper, we propose a Diffusion Posterior Sampling (DPS) strategy for solving general PDE-based inverse problems, combined with a subspace diffusion technique for sampling acceleration. This novel approach not only leverages the adjoint-state method with diffusion generative priors to solve PDE-based inverse problems but also enables PDE-constrained conditional sampling across subspace sequences, effectively reducing acquisition times. This approach effectively handles the relationship between the diffusion sampling part and the PDE inverse problem-solving part, designing corresponding subspace-based techniques to accelerate both processes simultaneously. We should note that this framework could be generalized to other PDE-based problems once the associated adjoint equations are explicitly known. Extensive experiments conducted on simulated speed field samples from the Marmousi and the KIT4 datasets, along with different geometries, show that the proposed method achieves significant improvements in travel-time imaging quality compared to existing methodologies. These findings highlight the potential of our proposed subspace DPS approach to advance the field of travel-time tomography, providing more accurate and high-quality speed field reconstructions.
\clearpage

\bibliographystyle{unsrt}
\bibliography{references}

\clearpage
\appendix
\section{Proofs} \label{sec:proof}

\bound* 
\noindent \textit{Proof.} Denote the rescaled speed field 
$$
\overline{\X}^k_0(u,v) := 2^k \X_0^k(2^{-k}u, 2^{-k}v)
$$
$$
\overline{\Y}^k_n(u,v) = \Y_n^k(2^{-k}u, 2^{-k}v)
$$
for $[u,v]\in [0,1]^2$, which also satisfy the Eikonal equation $\overline{\X}^k_0(u,v) \| \nabla \overline{\Y}^k_n(u, v)\|_2 = 1$ with the same domain of $\X_0(u,v)$ and $\Y_n(u,v)$. We continue to apply the Fourier transform to $\overline{\X}^k_0(u,v)$  to formulate
\begin{equation}\label{mid_result1}
    \begin{aligned}
        \mathcal{F}[\overline{\X}^k_0](\xi, \eta) &= \int_{\mathbb{R}^2} 2^k e^{-2\pi i (u\xi + v\eta)} \X_0^k(2^{-k}u, 2^{-k}v) \, du dv \\
        &= \int_{\mathbb{R}^2} 2^k e^{-2\pi i (u\xi + v\eta)} \Phi(\X^{k-1}_0)(2^{-k}u, 2^{-k}v) du dv\\
        &= \int_{\mathbb{R}^2} 2^{k-1} e^{-2\pi i (u\xi + v\eta)} du dv \left( \int_{\mathbb{R}^2} \varphi(\tilde{u}, \tilde{v})\X^{k-1}_0(2^{-k+1}u - \tilde{u}, 2^{-k+1}v - \tilde{v}) d\tilde{u} d\tilde{v}\right).\\
    \end{aligned}
\end{equation}
We substitute $(\tilde{u},\tilde{v})$ by $(2^{-k+1}\tilde{u},2^{-k+1}\tilde{v})$, and exchange the integration order to get
\begin{equation}
    \begin{aligned}
        \mathcal{F}[\overline{\X}^k_0](\xi, \eta) & = \left(\int_{\mathbb{R}^2}  \varphi(2^{-k+1}\tilde{u},2^{-k+1}\tilde{v}) e^{-2\pi i (\tilde{u}\xi + \tilde{v}\eta)} d (2^{-k+1}\tilde{u}) d(2^{-k+1}\tilde{v}) \right) \\
        & \left(\int_{\mathbb{R}^2}  \overline{\X}^{k-1}_0(u - \tilde{u}, v - \tilde{v})e^{-2\pi i ((u - \tilde{u})\xi + (v - \tilde{v})\eta)}  du dv\right) \\
        &= \hat{\varphi}(2^{k-1}\xi, 2^{k-1}\eta) \mathcal{F}[\overline{\X}^{k-1}_0](\xi, \eta) \\
    \end{aligned}
\end{equation}
where $\hat{\varphi}$ denotes the Fourier transformation of $\varphi$. Therefore, we can conclude that
\begin{equation}
    \begin{aligned}
        \mathcal{F}[\overline{\X}^k_0](\xi, \eta) & = \left(\prod\limits_{l=0}^{k-1}\hat{\varphi}(2^{l}\xi, 2^{l}\eta)\right) \mathcal{F}[\X_0](\xi, \eta) 
    \end{aligned}
\end{equation}
According to the Parseval's Identity in Sobolev space $H_\mu^1\left(\Omega\right)$, we have 
\begin{equation}
      \begin{aligned}
          \|\overline{\X}^k_0 - \X_0\|_{H_\mu^1\left(\Omega\right)}^2 &= \int_{\mathbb{R}^n} \left(1 + \mu(|\xi|^2 + |\eta|^2)\right) \left|\mathcal{F}[\overline{\X}^k_0] (\xi, \eta) - \mathcal{F}[\X_0](\xi, \eta)\right|^2 \,d\xi d\eta \\
          &=  \int_{\mathbb{R}^n} \left(1 + \mu(|\xi|^2 + |\eta|^2)\right) \left(\prod\limits_{l=0}^{k-1}\hat{\varphi}(2^{l}\xi, 2^{l}\eta) - 1\right)^2\left|\mathcal{F}[\X_0](\xi, \eta)\right|^2 \,d\xi d\eta \\
          & \le \int_{\mathbb{R}^n} \left(1 + \mu(|\xi|^2 + |\eta|^2)\right) \varepsilon^2 \left(1 + \mu(|\xi|^2 + |\eta|^2)\right)^r \left|\mathcal{F}[\X_0](\xi, \eta)\right|^2 \,d\xi d\eta \\
          & = \varepsilon^2 \int_{\mathbb{R}^n} \left(1 + \mu(|\xi|^2 + |\eta|^2)\right)^{r+1} \left|\mathcal{F}[\X_0](\xi, \eta)\right|^2 \,d\xi d\eta  \\
          & = \varepsilon^2 \|\X_0\|_{H_\mu^{1+r}\left(\Omega\right)}^2
      \end{aligned}
\end{equation}
By the definition of $\overline{\X}^k_0$, we have $\mathcal{J}_k(\X^k_0) = \mathcal{J}(\overline{\X}^k_0)$. If we regard $\overline{\X}^k_0$ as the perturbed speed field from $\X_0$, the change between $\mathcal{J}(\overline{\X}^k_0)$ and $\mathcal{J}(\X_0)$ is formulated as:
\begin{equation}
      \begin{aligned}
          \left|\mathcal{J}_k(\X^k_0)- \mathcal{J}(\X_0)\right| &= \left|\mathcal{J}(\overline{\X}^k_0)- \mathcal{J}(\X_0)\right| \\
          & = \partial \mathcal{J}(\X_0)(\overline{\X}^k_0 - \X_0)  +  O(\|\overline{\X}^k_0 - \X_0\|^2)  \\
          & =  \partial \mathcal{J}(\X_0)(\overline{\X}^k_0 - \X_0)  +  O( \varepsilon^2)
      \end{aligned}
\end{equation}
Apply the regularized adjoint-state method introduced in Section \ref{Adjoint_Method}, we formulate the first as
\begin{equation}
      \begin{aligned}
          \partial \mathcal{J}(\X_0)(\overline{\X}^k_0 - \X_0) &\le \|\partial \mathcal{J}(\X_0)\|_{H_\mu^1\left(\Omega\right)} \|\overline{\X}^k_0 - \X_0\|_{H_\mu^1\left(\Omega\right)} \\ 
          & \le \varepsilon \|\tilde{\X}_0\|^2_{H_\mu^1\left(\Omega\right)} \|\X_0\|_{H_\mu^{1+r}\left(\Omega\right)},\\
      \end{aligned}
\end{equation}
which finishes the proof. $\hfill\blacksquare$ 
    
\kernel*  
\noindent \textit{Proof.} Substitute the kernel function $\varphi_h(u,v)$ into 
\begin{equation}
    \begin{aligned}
        \X^{k+1}_0(u,v) &= \Phi(\X^{k}_0)(u,v) \\
        & = \frac{1}{2} \int_{\mathbb{R}^2}\varphi_h(\tilde{u},\tilde{v})\X^{k}_0(2u - \tilde{u}, 2v -\tilde{v}) d\tilde{u} d\tilde{v} \\
        & =   \frac{1}{2} \int_{\mathbb{R}^2} \frac{1}{4}\sum\limits_{(a,b)\in \{0,1\}^2} \delta(\tilde{u}+ah,\tilde{v}+bh) \X^{k}_0(2u - \tilde{u}, 2v -\tilde{v}) d\tilde{u} d\tilde{v} \\ 
        & = \frac{1}{2}\left(\sum\limits_{(a,b)\in \{0,1\}^2} \frac{\X^{k}_0(2u +ah, 2v +bh)}{4}\right).
    \end{aligned}
\end{equation}
Thus, by the definition of $\x^{k+1}_0$, we have
\begin{equation}
    \begin{aligned}
        \x_0^{k+1}[i,j] &= \X^k_0(i h, j h) \\
        &= \frac{1}{2}\left(\sum\limits_{(a,b)\in \{0,1\}^2} \frac{\X^{k}_0\left((2i+a)h, (2j+b)h\right)}{4}\right) \\ 
        &= \frac{1}{2}\left(\sum\limits_{(a,b)\in \{0,1\}^2} \frac{\x_0^{k}[2i+a, 2j+b]}{4}\right). 
    \end{aligned}
\end{equation}
$\hfill\blacksquare$ 

\bounddiscrete*  
\noindent \textit{Proof.} Consider the Fourier transform of the kernel function $\varphi_h(u,v)$ as
\begin{equation}
    \begin{aligned}
        \hat{\varphi}_h(\xi,\eta) &= \int_{\mathbb{R}^2} \varphi_h(u,v) e^{-2\pi i (u\xi + v\eta)}  \, du dv  \\ 
        &= \int_{\mathbb{R}^2} \frac{1}{4}\sum\limits_{(a,b)\in \{0,1\}^2} \delta(u+ah, v+bh) e^{-2\pi i (u\xi + v\eta)}  \, du dv \\
        &= \frac{1}{4}\left(1 + e^{2\pi h\xi i} + + e^{2\pi h\eta i} + + e^{2\pi h(\xi + \eta) i} \right) \\
        & = \frac{1}{2}\left[\cos{h(\xi + \eta)\pi} + \cos{h(\xi - \eta)\pi}\right] e^{\pi h(\xi + \eta)i} \\
        & =  \cos{(h\pi\xi)}\cos{(h\pi\eta)} e^{h\pi(\xi + \eta)i}, 
    \end{aligned}
\end{equation}
which indicates that $|\hat{\varphi}_h(\xi,\eta)|<1$. Then, we have the following inequality:
\begin{equation}
    \begin{aligned}
        \left|\hat{\varphi}_h(\xi,\eta) - 1\right| & = \left| \cos{h\pi(\xi+\eta)} e^{h\pi(\xi + \eta)i} - 1 + \sin{(h\pi\xi)}\sin{(h\pi\eta)} e^{h\pi(\xi + \eta)i}\right| \\ 
        & \le  \left| \cos{h\pi(\xi+\eta)} e^{h\pi(\xi + \eta)i} - 1 \right|  + \left| \sin{(h\pi\xi)}\sin{(h\pi\eta)} e^{h\pi(\xi + \eta)i}\right|  \\
        & \le |\sin{h\pi(\xi+\eta)}| + \frac{1}{2}|\sin{(h\pi\xi)}| |\sin{(h\pi\eta)}| + \frac{1}{2}|\sin{(h\pi\xi)}| |\sin{(h\pi\eta)}|\\
        & \le h\pi |\xi+\eta| + \frac{1}{2} h\pi |\xi| + \frac{1}{2} h\pi |\eta| \\
        & \le 3  h\pi \frac{|\xi| + |\eta|}{2} \\
        & \le \frac{3h\pi}{\sqrt{2}} (|\xi|^2 + |\eta|^2)^{\frac{1}{2}} 
    \end{aligned}
\end{equation}
According to this inequality, we verify the condition of the Theorem \ref{bound} by  
\begin{equation}\label{qwe}
    \begin{aligned}
        \left|\prod\limits_{l=0}^{k-1}\hat{\varphi}_h(2^{l}\xi, 2^{l}\eta) - 1\right| & = \left| \sum\limits_{s = 1}^{k-1}\left(\prod\limits_{l=0}^{s}\hat{\varphi}_h(2^{l}\xi, 2^{l}\eta) - \prod\limits_{l=0}^{s-1}\hat{\varphi}_h(2^{l}\xi, 2^{l}\eta) \right) + \left(\hat{\varphi}_h(\xi,\eta) - 1\right) \right|  \\ 
        & \le \sum\limits_{s = 1}^{k-1} \left|\hat{\varphi}_h(2^{s}\xi, 2^{s}\eta) - 1\right| \prod\limits_{l=0}^{s-1}\left|\hat{\varphi}_h(2^{l}\xi, 2^{l}\eta)\right| + \left|\hat{\varphi}_h(\xi,\eta) - 1\right| \\
        & \le \sum\limits_{s = 0}^{k-1}\left|\hat{\varphi}_h(2^{s}\xi, 2^{s}\eta) - 1\right| \\
        & \le \frac{3kh\pi}{\sqrt{2}} \sum\limits_{s = 0}^{k-1} (|2^{s}\xi|^2 + |2^{s}\eta|^2)^{\frac{1}{2}} \\
        & = \frac{3\pi k(2^k-1)}{\sqrt{2}} h (|\xi|^2 + |\eta|^2)^{\frac{1}{2}} \\
        & \le C_{k,\mu} h(1 + \mu(|\xi|^2 + |\eta|^2))^{\frac{1}{2}}
    \end{aligned}
\end{equation}
where $C_{k,\mu} = \frac{3\pi k(2^k-1)}{\sqrt{2\mu}}$ and $h$ is chosen such that $h < \frac{\varepsilon}{ C_{k,\mu}}$. So, Theorem \ref{bound} concludes that 
\begin{equation}
     \left|\mathcal{J}_k(\X^k_0)- \mathcal{J}(\X_0)\right| \le C_{k,\mu} h \|\tilde{\X}_0\|^2_{H_\mu^1\left(\Omega\right)} \|\X_0\|_{H_\mu^{2}\left(\Omega\right)}
\end{equation} 
$\hfill\blacksquare$ 

\bounded* 
\noindent \textit{Proof.} Theorem 2.6 in \cite{deckelnick2011numerical} concludes that for the bounded Lipschitz continuous speed field $\X_0^k$, the continuous viscosity solution $\Y_n^k$ in \eqref{eikonal_n} and the numerical solution $\y_n^k$ in \eqref{local_solver_section} have an $O(h^{\frac{1}{2}})$ error bound for enough small $h$ as follows:
\begin{equation}
    \max\limits_{(ih,jh)\in \mathbb{Z}_h^2}|\y_n^k[i,j] - \Y_n^k[ih,jh]| < C \sqrt{h}
\end{equation}
where $C$ depends on the bounds and Lipschitz continuity of $\X_0^k$, also the domain $\Omega$ and the regularity of the boundary $\partial \Omega$. Thus, given the definition of $\mathcal{J}^h_k(\x^k_0)$ in \eqref{misfit_approx_discrete} and $\mathcal{J}_k(\X^k_0)$ in \eqref{surrogate_misfit_k}, we have
\begin{equation}
    \begin{aligned}
        |\mathcal{J}^h_k(\x^k_0) - \mathcal{J}_k(\X^k_0)| & = \left|\frac{1}{2} \sum\limits_{m=1}^M\sum\limits_{n=1}^N \left(\y^{\delta}_{m,n} - \y^k_n[i^{r,k}_m,j^{r,k}_m]\right)^2 - \left(\y^{\delta}_{m,n} - \Y_n^k\left(\frac{u_m^r}{2^k}, \frac{v_m^r}{2^k}\right)\right)^2 \right|\\
        & \le \sum\limits_{m=1}^M\sum\limits_{n=1}^N \left|\Y_n^k\left(\frac{u_m^r}{2^k}, \frac{v_m^r}{2^k}\right) - \y^k_n[i^{r,k}_m,j^{r,k}_m]\right| \left|\y^{\delta}_{m,n} - \Y_n^k\left(\frac{u_m^r}{2^k}, \frac{v_m^r}{2^k}\right)\right| \\
        & + \frac{1}{2}  \sum\limits_{m=1}^M\sum\limits_{n=1}^N \left|\Y_n^k\left(\frac{u_m^r}{2^k}, \frac{v_m^r}{2^k}\right) - \y^k_n[i^{r,k}_m,j^{r,k}_m]\right|^2 \\
        & \le C\sqrt{h} \sum\limits_{m=1}^M\sum\limits_{n=1}^N \left|\y^{\delta}_{m,n} - \Y_n^k\left(\frac{u_m^r}{2^k}, \frac{v_m^r}{2^k}\right)\right| + \frac{1}{2}MNC^2h \\
        & \le  \sqrt{2MN}C\sqrt{h} \mathcal{J}_k(\X^k_0)^{\frac{1}{2}} + \frac{1}{2}MNC^2h \\
        & = C_1 \sqrt{h} + O(h)
    \end{aligned}
\end{equation}
where we assume that $(i^{r,k}_mh,j^{r,k}_mh) = ( \frac{u_m^r}{2^k}, \frac{v_m^r}{2^k})$ for simplicity, and $C_1$ depends on $\X^k_0$ and $\y^{\delta}$. Thus, we can conclude that
\begin{equation}
    \begin{aligned}
        |\mathcal{J}^h_k(\x^k_0) - \mathcal{J}(\X_0)| & = |\mathcal{J}^h_k(\x^k_0) - \mathcal{J}_k(\X^k_0) + \mathcal{J}_k(\X^k_0) - \mathcal{J}(\X_0)|\\
        & \le |\mathcal{J}^h_k(\x^k_0) - \mathcal{J}_k(\X^k_0)| + |\mathcal{J}_k(\X^k_0) - \mathcal{J}(\X_0)| \\
        & \le C_{k,\mu} h \|\tilde{\X}_0\|^2_{H_\mu^1\left(\Omega\right)} \|\X_0\|_{H_\mu^{2}\left(\Omega\right)} + C_1 \sqrt{h} + O(h) 
    \end{aligned}
\end{equation}
$\hfill\blacksquare$

\end{document}

%% file: math_commands.tex
\usepackage{amsfonts,bm}

\def\eqref#1{equation~\ref{#1}}
\def\1{\bm{1}}

\DeclareMathAlphabet{\mathsfit}{\encodingdefault}{\sfdefault}{m}{sl}
\SetMathAlphabet{\mathsfit}{bold}{\encodingdefault}{\sfdefault}{bx}{n}

\newcommand{\R}{\mathbb{R}}

\newcommand{\X}{{\mathbf X}}

\newcommand{\x}{{\boldsymbol x}}

\newcommand{\s}{{\boldsymbol s}}
\newcommand{\y}{{\boldsymbol y}}

\newcommand{\Y}{{\mathbf Y}}

\newcommand{\n}{{\boldsymbol n}}

\newcommand{\w}{{\boldsymbol w}}

\newcommand{\Rd}{{\mathbb R}}

\newcommand{\Ed}{{\mathbb E}}

\newcommand{\Ac}{{\mathcal A}}

\newcommand{\Nc}{{\mathcal N}}

%% file: packages.tex
\usepackage[utf8]{inputenc} % allow utf-8 input
\usepackage[T1]{fontenc}    % use 8-bit T1 fonts
\usepackage{hyperref}       % hyperlinks
\usepackage{url}            % simple URL typesetting
\usepackage{booktabs}       % professional-quality tables
\usepackage{amsfonts}       % blackboard math symbols
\usepackage{nicefrac}       % compact symbols for 1/2, etc.
\usepackage{microtype}      % microtypography
\usepackage{xcolor}       % colors
\usepackage{xkcdcolors}
\usepackage{times}
\usepackage{epsfig}
\usepackage{graphicx}
\usepackage{placeins}
\usepackage{multirow}
\usepackage[skip=5pt]{caption}
\usepackage{rotating}
\usepackage{cancel}
\usepackage{setspace}
\usepackage{eso-pic}

\usepackage{bm}
\usepackage{bibunits} % for separate ref list in appendix

\usepackage{amssymb,amsthm}
\usepackage{hyperref}
\usepackage{amsmath}
\usepackage{graphicx}
\usepackage{wrapfig,lipsum,booktabs}

\usepackage{epsfig}
\usepackage{tikz}
\usetikzlibrary{spy}
\usepackage{algpseudocode}
\usepackage{mathrsfs}

\usepackage{nicefrac}       % compact symbols for 1/2, etc.
\usepackage{booktabs}       % professional-quality tables

\usepackage{thmtools,thm-restate}
\usepackage{cleveref}

\usepackage{subcaption}        % subfigure